\documentclass[11pt]{article}

\usepackage[a4paper,margin=1in]{geometry}
\usepackage{amsmath,amssymb,amsthm,mathtools}
\usepackage{mathrsfs}
\usepackage{enumitem}
\usepackage{microtype}
\usepackage{booktabs}
\usepackage{array}
\usepackage{aliascnt}
\usepackage[hidelinks]{hyperref}
\usepackage{cleveref}
\hypersetup{
  pdftitle={High-Degree Restricted Cohomology for K(3;1)(1) in Characteristics Five and Seven},
  pdfauthor={Chao Ma},
  pdfcreator={Chao Ma},
  pdfsubject={High-degree restricted cohomology of the positive graded subalgebra of the rank-one contact Lie algebra, with parabolically induced coefficients, in characteristics five and seven},
  pdfkeywords={restricted cohomology, contact Lie algebra, K(3;1), positive graded subalgebra, parabolically induced module, central character, restricted Heisenberg algebra, PBW duality, Tor, Chevalley-Eilenberg complex, May resolution, Lyndon-Hochschild-Serre spectral sequence, modular representation theory}
}

\theoremstyle{plain}
\newtheorem{theorem}{Theorem}[section]
\newaliascnt{proposition}{theorem}
\newtheorem{proposition}[proposition]{Proposition}
\aliascntresetthe{proposition}
\newaliascnt{lemma}{theorem}
\newtheorem{lemma}[lemma]{Lemma}
\aliascntresetthe{lemma}
\newaliascnt{corollary}{theorem}
\newtheorem{corollary}[corollary]{Corollary}
\aliascntresetthe{corollary}
\theoremstyle{definition}
\newaliascnt{definition}{theorem}

\aliascntresetthe{definition}
\theoremstyle{remark}
\newaliascnt{remark}{theorem}

\aliascntresetthe{remark}
\numberwithin{equation}{section}

\newcommand{\kk}{k}
\newcommand{\F}{\mathbb F}
\newcommand{\Kalg}{\mathfrak K}
\newcommand{\slie}{\mathfrak{sl}_2}
\newcommand{\Ext}{\operatorname{Ext}}
\newcommand{\Tor}{\operatorname{Tor}}
\newcommand{\Hom}{\operatorname{Hom}}
\newcommand{\im}{\operatorname{im}}
\newcommand{\coker}{\operatorname{coker}}
\newcommand{\gr}{\operatorname{gr}}
\newcommand{\uenv}{\mathrm u}

\newcommand{\spb}{\mathrm{sp}}
\newcommand{\Dlt}{\Delta}

\title{High-Degree Restricted Cohomology for $K(3;\mathbf1)^{(1)}$ in Characteristics Five and Seven}
\author{Chao Ma\thanks{Independent Researcher; email: \texttt{raymond.ma@me.com}}}
\date{}

\begin{document}
\maketitle

\begin{abstract}
We determine the high-degree restricted first cohomology of the positive
subalgebra of the rank-one contact algebra in characteristics five and seven,
with parabolically induced coefficients. In degrees from twice the
characteristic through four times the characteristic minus three, each
characteristic has six nonzero homogeneous components, of total dimension thirteen.
The constituents include Frobenius twists of the natural horizontal module,
represented by Cartier cocycles. PBW duality converts the calculation into
bounded-degree homology, where the coefficient action couples the natural
and cubic relation coordinates. The resulting modules determine the
high-degree part of the coefficient cohomology used in the
Lyndon--Hochschild--Serre term $E_2^{0,1}$ for extensions of induced contact modules.
In particular, extensions from the Steinberg induced modules vanish in these degrees.
\end{abstract}

\section{Introduction}

Benkart and Feldvoss ask for the extensions between simple modules of the
restricted simple Lie algebras of Cartan type \cite[Problem~4]{BenkartFeldvoss}.
Hochschild introduced the restricted cohomology of a restricted Lie algebra
\cite{Hochschild}, and May constructed the free resolution of the trivial module
that we use \cite{May}.  Holmes determined the dimensions of the simple restricted
modules of the restricted contact algebras \cite{HolmesContact}.  Holmes and Zhang
\cite{HolmesZhang} studied the simple restricted modules of the Cartan-type algebras
through modules induced from the nonnegative part, and Zhang \cite{ZhangCartan}
constructed in the same way the simple modules with nonsingular $p$-characters of height
greater than one.  Lin and
Nakano related extensions between such modules to the cohomology of the positive graded
subalgebra \cite{LinNakano}.
For nonzero $p$-characters, we studied extensions through Peirce corners
in \cite{MaPeirce}. Here we work at zero $p$-character and compute the
coefficient cohomology entering the parabolic spectral sequence.

Let $\kk$ be a field of characteristic $p\in\{5,7\}$.
Let $\Kalg=K(3;\mathbf1)^{(1)}$, let $\mathfrak n=\Kalg_{>0}$, and let $X$ be a restricted
$\Kalg$-module.  For the induced modules $M(d,\zeta)$ defined below, PBW freeness and the
Eckmann--Shapiro isomorphism give
\[
 \Ext^1_{\uenv_0(\Kalg)}\bigl(M(d,\zeta),X\bigr)
 \cong
 \Ext^1_{\uenv_0(\Kalg_{\ge0})}\bigl(U_{d,\zeta},X\bigr),
\]
and for $0\to\mathfrak n\to\Kalg_{\ge0}\to\Kalg_0\to0$ the restricted
Lyndon--Hochschild--Serre spectral sequence has
\begin{equation}\label{eq:LHS}
 E_2^{i,j}
 =\Ext^i_{\uenv_0(\Kalg_0)}
 \bigl(U_{d,\zeta},H^j_{\mathrm{res}}(\mathfrak n,X)\bigr)
 \Longrightarrow
 \Ext^{i+j}_{\uenv_0(\Kalg_{\ge0})}(U_{d,\zeta},X).
\end{equation}
For $X$ an induced module, Theorem~\ref{thm:main} below computes
$H^1_{\mathrm{res}}(\mathfrak n,X)$ in the degrees $2p\le q\le4p-3$, and hence the
corresponding part of $E_2^{0,1}$, in characteristics five and seven.

The degree-zero subalgebra is
\[
 \Kalg_0\cong\slie\oplus\kk z.
\]
Let $V$ be the natural two-dimensional $SL_2$-module.  For $0\le d\le p-1$, let $L(d)$
denote the restricted simple $SL_2$-module of highest weight $d$; for a Weyl module of
possibly non-$p$-restricted highest weight we write $\Dlt(d)$.  For $\zeta\in\F_p$, let
$\kk_\zeta$ be the one-dimensional $\kk z$-module on which $z$ acts by $\zeta$, and put
\[
 U_{d,\zeta}=L(d)\boxtimes\kk_\zeta,
 \qquad
 M(d,\zeta)=\uenv_0(\Kalg)\otimes_{\uenv_0(\Kalg_{\ge0})}U_{d,\zeta},
\]
where $\Kalg_{>0}$ acts trivially on $U_{d,\zeta}$ and $\uenv_0$ denotes the restricted
enveloping algebra.  We write
\[
 H^j_{\mathrm{res}}(\mathfrak n,X)
 :=\Ext^j_{\uenv_0(\mathfrak n)}(\kk,X),
\]
with a second subscript for the homogeneous degree.

We consider
\[
 q=q_0+2p,
 \qquad
 0\le q_0\le2p-3.
\]
The negative part of $\Kalg$ is a restricted Heisenberg algebra with PBW top depth
\[
 D_{\mathrm{PBW}}=4(p-1).
\]
For a homogeneous chain in the dual Tor complex
$m\otimes x_1\wedge\cdots\wedge x_i$ we use the internal degree
\[
 \deg_{\mathrm{int}}(m\otimes x_1\wedge\cdots\wedge x_i)
 =\deg_{\mathrm{PBW}}m-\sum_j\deg_{\Kalg}x_j.
\]
Thus a degree-$q$ cochain is dual to the single internal degree
\[
 r=D_{\mathrm{PBW}}-q=2p-4-q_0.
\]
The central label changes to
\[
 \zeta^\dagger=-\zeta-4.
\]
For a parabolically induced module $M(\mu^\dagger)$, write $M(\mu^\dagger)_{\mathrm{rt}}$ for the right contragredient module obtained through the antipode.  We grade
\(\Tor^{\uenv_0(\mathfrak n)}_1(M(\mu^\dagger)_{\mathrm{rt}},\kk)\)
by the internal degree $r$.  Thus the homological degree is always $1$, whereas $r$ may
be equal to $-1$.  If a PBW coefficient is $x^ay^bt^c$, we call $c$ its $t$-exponent; $c=0$ is the
$t$-free face and $c\ge1$ the $t$-positive part.  On the $t$-free face set $h=r+1$.
Integer central-character labels are read in $\F_p$.  PBW duality gives a bounded internal range. We compute the coupled
ordinary coefficient relations in that range and then impose the restricted
relations. The nonzero classes are detected by tensor-module extensions and
closed differential forms.

\begin{theorem}[High-degree restricted cohomology]\label{thm:main}
Let $p\in\{5,7\}$, $0\le d<p$, and $\zeta\in\mathbb F_p$.
For $q=q_0+2p$ with $0\le q_0\le2p-3$, put
\[
 h=4p-3-q,\qquad r=h-1,\qquad
 \zeta^\dagger=-\zeta-4,\qquad \lambda=\zeta+4=-\zeta^\dagger.
\]
There is a graded duality
\begin{equation}\label{eq:main-duality}
 H^1_{\mathrm{res},q}(\mathfrak n,M(d,\zeta))^\vee
 \cong
 \Tor^{\uenv_0(\mathfrak n)}_{1,h-1}
       (M(d,\zeta^\dagger)_{\mathrm{rt}},\kk).
\end{equation}
The nonzero cohomology groups, as horizontal $SL_2$-modules, are
\begin{equation}\label{eq:corrected-main-table}
\begin{array}{c|c|c|c}
 h&d&\lambda&H^1_{\mathrm{res},q}(\mathfrak n,M(d,\zeta))\\ \hline
 0&0&0&L(1)\oplus L(3)\\
 0&1&1&L(0)\\
 0&3&1&L(0)\\
 1&0&0&L(0)\\
 1&p-2&0&V^{(1)}\\
 p+1&0&0&V^{(1)}.
\end{array}
\end{equation}
All other groups in this range vanish.
Here $V^{(1)}$ is the Frobenius twist of the natural two-dimensional
$SL_2$-module; its restricted $\slie$-action is trivial.
In particular $H^1_{\mathrm{res},2p}(\mathfrak n,M(d,\zeta))=0$
for every inducing label $(d,\zeta)$.
\end{theorem}

\begin{corollary}\label{cor:congruence-row}
Under the congruence
\begin{equation}\label{eq:central-character-main}
 \zeta^\dagger\equiv h\pmod p,
 \qquad\text{equivalently}\qquad\zeta\equiv q_0-1\pmod p,
\end{equation}
the only nonzero group in Theorem~\ref{thm:main} occurs at
$h=0$, $d=0$, $\zeta^\dagger=0$, and is $L(1)\oplus L(3)$.
\end{corollary}
\begin{proof}
Substitute the six labels of \eqref{eq:corrected-main-table} into
\eqref{eq:central-character-main}.
\end{proof}

\begin{theorem}[Cartier classes]\label{thm:qtop}
At $q=3p-4$, $d=0$, $\zeta=p-4$, the classes arising from the closed forms
\[
 x^{p-1}\,dx,\qquad y^{p-1}\,dy
\]
span a copy of $V^{(1)}$ in restricted first cohomology.
Under \eqref{eq:main-duality} these cocycles pair nondegenerately with
\[
 C_3\otimes x^2y^{p-1}u_0,\qquad
 C_0\otimes x^{p-1}y^2u_0,
\]
and with the truncation cycles
\[
 b_+\otimes x^{p-1}tu_0,\qquad
 b_-\otimes y^{p-1}tu_0.
\]
\end{theorem}

We grade $\Ext^1$ by the cochain convention above: a degree-$q$ cochain
on positive inputs of total contact degree $J$ takes values in PBW depth $q-J$.

\begin{corollary}[Steinberg inducing module]\label{cor:steinberg-source}
Let $p=5,7$, $\alpha,\zeta\in\mathbb F_p$, and $0\le d<p$.
For $2p\le q\le4p-3$,
\[
 \Ext^1_{\uenv_0(\Kalg),q}
       (M(p-1,\alpha),M(d,\zeta))=0.
\]
\end{corollary}
\begin{proof}
After extending scalars to an algebraic closure, the Steinberg module
$L(p-1)$ is projective over
$\uenv_0(\slie)$; see \cite[proof of Proposition~5, p.~3]{McNinchJordan}.
The algebra $\uenv_0(\kk z)=\kk[z]/(z^p-z)$ is semisimple, so
$U_{p-1,\alpha}$ is projective over $\uenv_0(\Kalg_0)$.
Thus \eqref{eq:LHS} is concentrated in column zero, and PBW freeness gives
\[
 \Ext^1_{\uenv_0(\Kalg),q}(M(p-1,\alpha),M(d,\zeta))
 \cong\Hom_{\uenv_0(\Kalg_0)}
       (U_{p-1,\alpha},H^1_{\mathrm{res},q}(\mathfrak n,M(d,\zeta))).
\]
By Theorem~\ref{thm:main}, the target has restricted horizontal types
$L(0)$, $L(1)$ and $L(3)$; $V^{(1)}$ restricts
to two trivial $\slie$-modules. Since $p-1\notin\{0,1,3\}$,
the Hom space is zero. The finite cochain complexes commute with
field extension, so vanishing descends to $\kk$.
\end{proof}

We establish PBW duality in Section~\ref{sec:pbw}, and then
describe the ordinary and restricted coefficient boundaries.
The first two internal rows are treated in Section~\ref{sec:first-two}.
For the other rows we use the highest representatives of
Section~\ref{sec:highest-data} and the common relation matrix in
Appendix~\ref{app:comparison}. Its few surviving families are
settled by the truncation obstructions of
Sections~\ref{sec:maxz-natural} and~\ref{sec:positive-fibres}.

\section{The contact algebra and PBW duality}\label{sec:contact}

For $\Kalg=K(3;\mathbf1)^{(1)}$ we use the truncated polynomial normalization of the divided-power realization; thus the monomials $x^ay^bz^c$, $0\le a,b,c\le p-1$, form the
coordinate basis and $x^p=y^p=z^p=0$ in this normalization.  This is the convention used
for the contact formulas below; see \cite{HolmesContact,ZhangCartan}.  For $p>3$ the
derived algebra $K(3;\mathbf1)^{(1)}$ is simple, hence centerless.  Put
\[
 \{f,g\}_{\spb}=\partial_xf\,\partial_yg-\partial_yf\,\partial_xg,
 \qquad
 \mathcal E=x\partial_x+y\partial_y.
\]
The contact bracket is
\begin{equation}\label{eq:contact-bracket}
 \{f,g\}_K
 =\{f,g\}_{\spb}
 +(2f-\mathcal Ef)\partial_zg
 +(\partial_zf)(\mathcal Eg-2g).
\end{equation}
The restricted map is characterized by
\[
 \operatorname{ad}(f^{[p]})=(\operatorname{ad}f)^p.
\]

The positive grading is the contact grading associated with $\deg x=\deg y=1$ and
$\deg z=2$.  In particular
\[
 \Kalg_1=C\oplus N,
 \qquad C=S^3V,
 \qquad N=zV,
\]
and we write $C_i=x^{3-i}y^i$ for $0\le i\le3$. The next component is
\[
 \Kalg_2=A\oplus B\oplus Q,
 \qquad A=S^4V,
 \qquad B=zS^2V,
 \qquad Q=\kk z^2.
\]
The negative part is the restricted Heisenberg algebra
\[
 \Kalg_-=\kk x\oplus\kk y\oplus\kk t,
 \qquad [x,y]=t,
\]
where $t$ is represented by the constant contact function.  Give PBW monomials the depth
\[
 \deg_{\mathrm{PBW}}x=\deg_{\mathrm{PBW}}y=1,
 \qquad
 \deg_{\mathrm{PBW}}t=2.
\]
Set
\[
 b_+=zx,
 \qquad b_-=zy,
 \qquad Q=z^2.
\]
Directly from \eqref{eq:contact-bracket},
\begin{equation}\label{eq:basic-contact}
 {[b_+,t]}=-2x,
 \qquad
 {[b_-,t]}=-2y,
 \qquad
 {[b_+,b_-]}=z^2,
\end{equation}
and
\begin{equation}\label{eq:basic-horizontal}
 {[b_+,x]}=-x^2,
 \quad [b_+,y]=z-xy,
 \quad [b_-,x]=-z-xy,
 \quad [b_-,y]=-y^2.
\end{equation}

\subsection{The restricted complex in low degree}\label{sec:restricted-complex}

Put $\mathfrak n=\Kalg_{>0}$ and $\mathcal A=\uenv_0(\mathfrak n)$.  We use the first three
terms of May's free resolution of the trivial restricted module \cite{May}.
In homological degrees $0,1,2$ they are
\begin{equation}\label{eq:may-low}
 P_0=\mathcal A,
 \qquad
 P_1=\mathcal A\otimes\mathfrak n,
 \qquad
 P_2=\mathcal A\otimes\bigl(\Lambda^2\mathfrak n\oplus\mathfrak n^{(1)}\bigr),
\end{equation}
where $\mathfrak n^{(1)}$ denotes the Frobenius-twisted copy carrying the restricted
$p$-power generators.  On a homogeneous basis element $x$, the last formula below is read
in $\mathfrak n^{(1)}$ and then extended $\kk$-linearly.  The differential is
\begin{align}
 d_1(a\otimes x)&=ax,\label{eq:may-d1}\\
 d_2(a\otimes x\wedge y)
 &=ax\otimes y-ay\otimes x-a\otimes[x,y],\label{eq:may-lie}\\
 d_2(a\otimes x^{(1)})
 &=ax^{p-1}\otimes x-a\otimes x^{[p]}.\label{eq:may-p}
\end{align}
May proves that these maps form the initial segment of a free $\mathcal A$-resolution; in
particular they are exact at $P_1$.  In the PBW presentation of $\mathcal A$, this says precisely
that the syzygies of the augmentation ideal are generated by the ordinary Lie relations
and the restricted relations $x^p-x^{[p]}$.  Thus
\eqref{eq:may-low}--\eqref{eq:may-p} compute $\Tor^{\mathcal A}_1(-,\kk)$.

For a right $\mathcal A$-module $X$, tensoring with \eqref{eq:may-low} gives
\begin{align}
 \partial_1(m\otimes x)&=mx,\label{eq:res-d1}\\
 \partial_2^{\mathrm{Lie}}(m\otimes x\wedge y)
 &=mx\otimes y-my\otimes x-m\otimes[x,y],\label{eq:res-lie}\\
 \partial_2^{[p]}(m\otimes x^{(1)})
 &=mx^{p-1}\otimes x-m\otimes x^{[p]}.\label{eq:res-p}
\end{align}
After transport through the antipode, \eqref{eq:res-lie} is the usual left
Chevalley--Eilenberg boundary, up to the global sign convention, while
\eqref{eq:res-p} is the restricted $p$-row.  Every ordinary differential used below is a
homogeneous component of \eqref{eq:res-lie}; every restricted differential in
Section~\ref{sec:restricted} is a homogeneous component of \eqref{eq:res-p}.  The ordinary three-boundaries used to change relation representatives are the
Chevalley--Eilenberg syzygies among the generators in \eqref{eq:res-lie}.

\subsection{PBW duality and the internal Tor grading}\label{sec:pbw}

Let
\[
 \Omega_-=x^{p-1}y^{p-1}t^{p-1}
\]
be the unique PBW-top monomial.  If $\operatorname{top}$ extracts its coefficient, set
\[
 \langle a,b\rangle_{\mathrm{PBW}}=\operatorname{top}(ab).
\]
Multiplication followed by projection to $\Omega_-$ pairs complementary PBW monomials
perfectly.  The induced contragredient duality on parabolically induced modules changes the central character by
$-4$.  Here $M^*$ denotes the usual left contragredient module, obtained from the
vector-space dual through the antipode.

\begin{proposition}[PBW top-coefficient duality]\label{prop:PBW-dual}
For $\mu=(d,\zeta)$ put
\[
 \mu^\dagger=(d,-\zeta-4).
\]
Then
\[
 M(\mu)^*\cong M(\mu^\dagger),
\]
and the depth-$s$ and depth-$(D_{\mathrm{PBW}}-s)$ pieces are perfectly paired.
\end{proposition}

\begin{proof}
The PBW top-coefficient form is perfect on $\uenv_0(\Kalg_-)$.
The element $z$ acts on $x,y,t$ with respective weights
$-1$, $-1$, and $-2$.
Hence their top PBW monomial has weight
\[ (p-1)(-1-1-2)=4\quad\text{in }\mathbb F_p .\]
The inducing line in the contragredient therefore has central
character $-\zeta-4$. The top line has trivial horizontal character,
and $L(d)$ is self-contragredient, giving the asserted inducing
module and the complementary depth pairing.
\end{proof}

If $M$ is a finite-dimensional left $\mathcal A$-module, its vector-space dual is a right
$\mathcal A$-module and finite-dimensional duality yields
\[
 \Ext^1_{\mathcal A}(\kk,M)^\vee\cong\Tor^{\mathcal A}_1(M^\vee,\kk).
\]
Via the antipode, the right dual of $M(\mu)$ is identified with the right-module form of
$M(\mu^\dagger)$.  Write $M(\mu^\dagger)_{\mathrm{rt}}$ for this right module.

\begin{theorem}[Graded Ext--Tor duality]\label{thm:pbw-duality}
Let
\[
 q=q_0+2p,
 \qquad 0\le q_0\le2p-3,
\]
and put
\[
 r=D_{\mathrm{PBW}}-q=2p-4-q_0.
\]
Then, degree by degree,
\begin{equation}\label{eq:pbw-duality}
 H^1_{\mathrm{res},q}(\mathfrak n,M(\mu))^\vee
 \cong
 \Tor^{\mathcal A}_{1,r}
 \bigl(M(\mu^\dagger)_{\mathrm{rt}},\kk\bigr)
 =
 \Tor^{\uenv_0(\mathfrak n)}_{1,r}
 \bigl(M(\mu^\dagger)_{\mathrm{rt}},\kk\bigr).
\end{equation}
Consequently, the cohomological range corresponds to the bounded internal interval
\[
 -1\le r\le2p-4.
\]
\end{theorem}

\begin{proof}
Let $P_\bullet\to\kk$ be a graded projective resolution over $\mathcal A$.
Finite-dimensional duality yields
\[
 \Hom_{\mathcal A}(P_i,M)^\vee\cong M^\vee\otimes_{\mathcal A}P_i,
\]
and the transpose of the cochain differential is the tensor differential.  By Section~\ref{sec:restricted-complex}, the first two tensor differentials are precisely
\eqref{eq:res-lie} and \eqref{eq:res-p}; transport through the antipode gives the left
convention used from this point onward.  A degree-$q$ cochain on a positive tensor of total
degree $J$ takes values in coefficient depth $q-J$.  After PBW duality this becomes depth
$D_{\mathrm{PBW}}-q+J$.  Subtracting the positive degree $J$ gives the constant internal degree
$D_{\mathrm{PBW}}-q=r$.
\end{proof}

\section{Relation modules and coefficient boundaries}\subsection{The \texorpdfstring{$t$}{t}-exponent filtration}\label{sec:koszul}

In the $\Kalg_1$-head part of the degree-one Tor term, a coefficient monomial has the
form
\[
 x^ay^bt^c\otimes u,
 \qquad
 a+b+2c=r+1.
\]
For a head in $\Kalg_j$ the coefficient depth is instead $r+j$.  By
Lemma~\ref{lem:degree-one-generation}, higher heads can be eliminated using the full
coefficient boundary \eqref{eq:left-boundary}; we retain its action terms in the relation
coordinates below.  On the $\Kalg_1$-head part, $r\le2p-4$ gives $c\le p-2$ and
$(a,b)\ne(p-1,p-1)$.  We use the filtration by the exponent of $t$ on these coordinates.

Let
\[
 R=\kk[x,y]/(x^p,y^p),
 \qquad R_m=R_m^{\mathrm{hom}}.
\]
For $c\ge1$, the degree $-1$ part of the differential for this filtration is obtained from
\eqref{eq:basic-contact}:
\begin{equation}\label{eq:t-drop-action}
 \sigma_t(b_+)\bigl(x^ay^bt^c\otimes u\bigr)
 =-2c\,x^{a+1}y^bt^{c-1}\otimes u,
\end{equation}
\begin{equation}\label{eq:t-drop-action-minus}
 \sigma_t(b_-)\bigl(x^ay^bt^c\otimes u\bigr)
 =-2c\,x^ay^{b+1}t^{c-1}\otimes u.
\end{equation}
The relation $[b_+,b_-]=z^2$ supplies the corresponding Koszul syzygy.

\begin{proposition}[Truncated Koszul complex]\label{prop:Koszul}
Fix $c\ge1$ and put
\[
 h=r+1-2c.
\]
Up to multiplication by a nonzero scalar, the leading $t$-lowering differential is
\[
 R_{h-1}\otimes U_{d,\zeta^\dagger}
 \xrightarrow{(-y,x)}
 R_h^{\oplus2}\otimes U_{d,\zeta^\dagger}
 \xrightarrow{(x,y)}
 R_{h+1}\otimes U_{d,\zeta^\dagger}.
\]
Its degree-one homology is zero unless $h=p-1$, in which case
\[
 \mathscr F_h
 =\kk[x^{p-1}e_+]\oplus\kk[y^{p-1}e_-].
\]
The only $t$-positive truncation degrees are therefore
\[
 p=5:\ q=11,
 \qquad
 p=7:\ q=15,17.
\]
\end{proposition}

\begin{proof}
The displayed complex is the homogeneous Koszul complex of $(x,y)$ on the complete
intersection $R=\kk[x,y]/(x^p,y^p)$.  Its first homology is generated by the truncation
relations $x^p=0$ and $y^p=0$, represented by $x^{p-1}e_+$ and $y^{p-1}e_-$.  Both lie in
horizontal degree $p-1$.  The degree-two Koszul generator has horizontal degree $2p-2$ and
cannot occur in a $t$-positive Tor slice.  The equation $h=p-1=r+1-2c$ gives
$r=p-2+2c$, and the permitted values of $c$ yield exactly the degrees displayed above.
\end{proof}

At $c=0$ the leading differential from the preceding $t$-degree is again the Koszul map.
The natural degree-one quotient fits into
\begin{equation}\label{eq:bottom-koszul}
 0\longrightarrow\delta_{h,p-1}\mathscr F_0
 \longrightarrow\mathscr N_h
 \longrightarrow R_{h+1}\longrightarrow0,
\end{equation}
where now $h=r+1$.

\subsection{Quadratic relations}\label{sec:quadratic}

The degree-two contact bracket
\[
 \beta:\Lambda^2\Kalg_1\longrightarrow\Kalg_2
\]
is surjective.  The three homogeneous summands behave as follows.

\begin{lemma}[Minimal quadratic relation module]\label{lem:quadratic-kernel}
Under
\[
 \Kalg_1=C\oplus N,
 \qquad
 \Kalg_2=A\oplus B\oplus Q,
\]
one has:
\begin{enumerate}[label=\textup{(\roman*)}]
\item $\Lambda^2N\to Q$ is an isomorphism;
\item $C\otimes N\to A\oplus B$ is an isomorphism;
\item $\Lambda^2C\to A$ is surjective with one-dimensional kernel.
\end{enumerate}
It follows that
\[
 \ker\beta\cong L(4)\oplus L(0).
\]
A generator of the scalar summand is
\begin{equation}\label{eq:r0}
 r_0=x^3\wedge y^3-3x^2y\wedge xy^2.
\end{equation}
\end{lemma}

\begin{proof}
Direct application of \eqref{eq:contact-bracket} to the indicated homogeneous summands gives (i)--(iii).  For the scalar relation one uses
\[
 \{x^3,y^3\}_{\spb}=9x^2y^2,
 \qquad
 \{x^2y,xy^2\}_{\spb}=3x^2y^2.
\]
The dimensions are $\dim\Lambda^2(C\oplus N)=15$ and $\dim\Kalg_2=9$, so the kernel has
dimension six; its two displayed $SL_2$ summands have dimensions five and one.
\end{proof}

By Lemma~\ref{lem:quadratic-kernel}, the restrictions
\[
 \beta_{CN}:C\otimes N\xrightarrow{\sim}A\oplus B,
 \qquad
 \beta_{NN}:\Lambda^2N\xrightarrow{\sim}Q
\]
are equivariant isomorphisms.  Define
\begin{equation}\label{eq:explicit-bracket-section}
 s
 =\beta_{CN}^{-1}\oplus\beta_{NN}^{-1}:
 (A\oplus B)\oplus Q\longrightarrow\Lambda^2\Kalg_1.
\end{equation}
Then $\beta s=\mathrm{id}_{\Kalg_2}$, so $s$ is an equivariant splitting.
In the ordinary differential \eqref{eq:res-lie}, the bracket
coordinate of $m\otimes s(y)$ is $\pm m\otimes y$ with unit coefficient.  The corresponding
summand is therefore contractible and may be removed.

The $L(4)$ relation produces the cubic transfer operator
\begin{equation}\label{eq:Ph}
 P_h:S^4V\otimes R_{h-1}\longrightarrow S^3V\otimes R_h,
 \qquad
 P_h(f\otimes m)
 =\partial_xf\otimes xm+\partial_yf\otimes ym.
\end{equation}

\begin{proposition}[Cubic transfer cokernel]\label{prop:cubic-transfer}
For every $p\ge5$,
\[
 \coker P_h\cong
 \begin{cases}
 L(3-h),&0\le h\le3,\\
 0,&h\ge4.
 \end{cases}
\]
\end{proposition}

\begin{proof}
Use $u,v$ for the symmetric-power variables and put
$T_{i,a}=u^{3-i}v^ix^ay^{h-a}$.
The nonzero coefficients satisfy
$0\le i\le3$ and $\max(0,h-p+1)\le a\le\min(h,p-1)$.
For a source monomial $S_{j,b}=u^{4-j}v^jx^by^{h-1-b}$,
\[
 P_h(S_{j,b})=(4-j)T_{j,b+1}+jT_{j-1,b}.
\]
Both terms have the same value of $a-i$; hence the matrix splits
into chains of at most four consecutive vertices. All internal
edge coefficients are units. If the left endpoint has $a>0$,
it is the sole nonzero term in the image of $S_{i,a-1}$; if the
right endpoint has $h-a>0$, it is the sole nonzero term in the
image of $S_{i+1,a}$. These sources satisfy the truncation bounds.
Starting at either such endpoint eliminates the whole chain.
A surviving chain must therefore start at $a=0$ and end at $a=h$,
which is possible only for $h\le3$, in exactly $4-h$ chains.

For $0\le h\le3$, the equivariant map
\[
 \pi_h=(\partial_u\partial_y-\partial_v\partial_x)^h
 :S^3V\otimes R_h\longrightarrow S^{3-h}V
\]
annihilates $P_h$ and is surjective, since
\[
 \pi_h(u^{3-i}v^iy^h)=(3-i)_h h!u^{3-h-i}v^i.
\]
This identifies the cokernel with $L(3-h)$.
For the matching highest vectors, writing $\omega=uy-vx$ gives
\[
 P_h(\omega^ju^{4-j}x^{h-1-j})
 =(4-j)\omega^ju^{3-j}x^{h-j}.
\]
Here $0\le j\le\min(3,h-1)$. The extra source
$\omega^4x^{h-5}$, when $h\ge5$, has image zero.
Thus the matching scalars are $4-j$; only $j\le h-1$ occurs
in the four low rows.
\end{proof}

The scalar relation \eqref{eq:r0} acts only on the cubic summand.  For $f\in S^3V$ and
$\varphi(x,y)\otimes u\in R_m\otimes L(d)$, its $t$-preserving action is
\begin{equation}\label{eq:cubic-action}
\begin{aligned}
 \mathcal D_f(\varphi\otimes u)
={}&-(\partial_x\varphi)\otimes\rho_d(\partial_yf)u
 +(\partial_y\varphi)\otimes\rho_d(\partial_xf)u\\
&+\frac12\Bigl(
 (\partial_y^2f)(\partial_x^2\varphi)
 -2(\partial_x\partial_yf)(\partial_x\partial_y\varphi)
 +(\partial_x^2f)(\partial_y^2\varphi)
 \Bigr)\otimes u.
\end{aligned}
\end{equation}
Projection of its boundary to $\coker P_h$ defines
\[
 \overline{\mathcal Q}_{h,d}:R_{h+1}\otimes L(d)
 \longrightarrow L(3-h)\otimes L(d).
\]
On the matching Weyl sections the scalars are:
\begin{equation}\label{eq:Q-scalar-table}
\begin{array}{c|c|c}
 h&j\to k&\text{scalar}\\ \hline
0&0\to1&-6d\\
0&1\to2&-6(d-1)\\ \hline
1&0\to0&-12(d+2)\\
1&1\to1&0\\
1&2\to2&12d\\ \hline
2&1\to0&24(d+3)\\
2&2\to1&24(d+2)\\ \hline
3&2\to0&0.
\end{array}
\end{equation}
Rows for which the relevant section is absent are omitted.

For the scalar relation, fix $\eta_j=(xY-yX)^jx^{h+1-j}X^{d-j}$ and
$\theta_k=(uY-vX)^ku^{3-h-k}X^{d-k}$, where $k=j+1-h$.
Apply $\pi_h$ to
\[
 v^3\mathcal D_{x^3}\eta_j-u^3\mathcal D_{y^3}\eta_j
 -3uv^2\mathcal D_{x^2y}\eta_j+3u^2v\mathcal D_{xy^2}\eta_j.
\]
Using \eqref{eq:cubic-action}, its coefficient on $\theta_k$ is
\[
\begin{array}{c|c}
 h&\lambda_{h,j}\\\hline
0&-6(d-j)\\
1&12\{j(d-j+3)-d-2\}\\
2&24\{jd-2j^2+5j-\binom j2(d-j+2)\}\\
3\ (j=2)&-216+(-144d+216)+(144d-216)+216=0.
\end{array}
\]
For $h=1$ the two contributions are
$18j(d-j+1)-18$ and
$-12d+24j-6-6j(d-j+1)$.
For $h=2$, putting $b=\binom j2$, the three contributions are
$72j-72b(d-j+2)$,
$12j(4d-8j+6)+48b(d-j+2)$, and
$-24jd+48j^2-24j$.
Taking the admissible $j$ gives every entry of
\eqref{eq:Q-scalar-table}.

\subsection{Relations in total degree three}\label{sec:cubic}

Let
\[
 \mathcal S_3(\mathfrak n)
 =H_2(\mathfrak n,\kk)_3
 =\frac{\ker(\Kalg_1\otimes\Kalg_2\to\Kalg_3)}
 {\im(\Lambda^3\Kalg_1\to\Kalg_1\otimes\Kalg_2)}.
\]
In characteristics five and seven, the Jacobi map is injective and the bracket map onto
$\Kalg_3$ is surjective.  The filtration by $z$-degree determines the associated graded
module.

\begin{proposition}[Degree-three relation module]\label{prop:cubic-relations}
Let $\chi_m=\operatorname{ch}\Dlt(m)$ denote the Weyl character of highest weight $m$.
Then
\[
 \operatorname{ch}\mathcal S_3=
 \begin{cases}
  \chi_7+\chi_5+\chi_3+2\chi_1,&p=7,\\
  \chi_7+\chi_5+\chi_3+2\chi_1+\operatorname{ch}V^{(1)},&p=5,
 \end{cases}
\]
and
\[
 \dim\mathcal S_3=
 \begin{cases}
 22,&p=7,\\
 24,&p=5.
 \end{cases}
\]
The $z$-degree filtration has associated-graded characters
\begin{equation}\label{eq:cubic-zdegree-table}
\begin{array}{c|c|c}
\text{$z$-degree}&p=7&p=5\\ \hline
3&\chi_1&\chi_1\\
2&\chi_3&\chi_3\\
1&\chi_5&\chi_5\\
0&\chi_7+\chi_1&\chi_7+\chi_1+\operatorname{ch}V^{(1)}.
\end{array}
\end{equation}
Moreover,
$\gr_2\mathcal S_3\cong L(3)$ and $\gr_3\mathcal S_3\cong L(1)$.
\end{proposition}

\begin{proof}
Write $\Kalg_1=C\oplus N$ and $\Kalg_2=A\oplus B\oplus Q$.  For
\[
 \mathcal J(a\wedge b\wedge c)
 =a\otimes[b,c]-b\otimes[a,c]+c\otimes[a,b],
\]
the source decomposes as
\[
 \Lambda^3\Kalg_1
 =\Lambda^3C\oplus(\Lambda^2C\otimes N)\oplus(C\otimes\Lambda^2N).
\]
The injectivity of $\mathcal J$ follows from three explicit projections used successively.  On
$C\otimes\Lambda^2N$, projection to $C\otimes Q$ is
$\mathrm{id}_C\otimes\beta_{NN}$ and is an isomorphism.  After this summand has been
removed, the composite of $\mathcal J|_{\Lambda^2C\otimes N}$ with projection onto $N\otimes A$ is
$\beta_{CC}\otimes\mathrm{id}_N$.  Its kernel is $\kk r_0\otimes N$.  On that kernel,
projection to $C\otimes B$ gives
\begin{align*}
 \pi_{C\otimes B}\mathcal J(r_0\wedge b_+)
 &=-3x^3\otimes zy^2
   +6x^2y\otimes zxy
   -3xy^2\otimes zx^2,
\end{align*}
which is nonzero in characteristics five and seven.  Since
$\kk r_0\otimes N\cong L(1)$ is simple, the restriction is injective.  Finally,
$\Lambda^3C\cong L(3)$ is simple and
\begin{align*}
 \mathcal J(x^3\wedge x^2y\wedge xy^2)
 &=3x^3\otimes x^2y^2
   -6x^2y\otimes x^3y
   +3xy^2\otimes x^4
\end{align*}
is nonzero.  Hence $J$ is injective on $\Lambda^3C$ as well, and therefore on all of
$\Lambda^3\Kalg_1$.

The bracket $\Kalg_1\otimes\Kalg_2\to\Kalg_3$ is surjective.  Its $z$-degree-three
component vanishes on $N\otimes Q$.  The lower target degrees are generated by
the nonzero equivariant brackets in $[N,B]$, $[C,B]+[N,A]$, and $[C,A]$; in characteristic
five the two pure fifth-power monomials do not occur in $\Kalg_3$.  The quotient of the
bracket kernel by the Jacobi image therefore has the four graded pieces in
\eqref{eq:cubic-zdegree-table}, whose sum gives the character and dimension claimed for
$\mathcal S_3$.  The $z$-degree-three piece is $N\otimes Q\cong L(1)$, while the $z$-degree-two quotient has
character $\chi_3$.  Since $L(3)$ is simple for $p=5,7$, the module identifications follow.
\end{proof}

\subsection{Coefficient boundaries}\label{sec:filtration}

Transport a right coefficient module through the antipode, so that
$a\cdot m=-ma$, and multiply chains of homological degree $j$ by $(-1)^j$.
We then use the left-coefficient convention
\begin{equation}\label{eq:left-boundary}
 \partial_2(m\otimes a\wedge b)
 =am\otimes b-bm\otimes a+m\otimes[a,b].
\end{equation}
The PBW action formulas below refer to this left action.

\begin{lemma}[Coefficient boundary identity]\label{thm:triangularity}
For $a,b,c\in\mathfrak n$,
\begin{equation}\label{eq:coefficient-boundary-change}
\begin{aligned}
 \partial_3(m\otimes a\wedge b\wedge c)={}&
 am\otimes b\wedge c-bm\otimes a\wedge c+cm\otimes a\wedge b\\
 &+m\otimes\bigl([a,b]\wedge c-[a,c]\wedge b+[b,c]\wedge a\bigr).
\end{aligned}
\end{equation}
If $a,b,c$ have total contact degree $J$, every bracket face has
that same total degree $J$.
\end{lemma}
\begin{proof}
Apply the ordinary Chevalley--Eilenberg differential in the convention
\eqref{eq:left-boundary}. The six terms give the displayed identity;
the Jacobi identity and the module identity give $\partial_2\partial_3=0$.
Contact degree is additive under each nonzero bracket.
\end{proof}

\section{Relations of total degree at least four}\label{sec:higher}

Let $F_JC_2(r)$ be the span of coefficient two-chains of internal degree $r$
whose two positive heads have total contact degree at most $J$.

For $J\ge4$ process ordinary relation degrees in increasing order and write
\[
 B_{<J}(r)=d_2\bigl(F_{J-1}C_2(r)\bigr).
\]
Lemma~\ref{thm:triangularity} allows us to change a degree-$J$ relation representative
by a complete Chevalley--Eilenberg three-boundary modulo $B_{<J}(r)$.

Put $\delta_t=\operatorname{ad}(t)$ and $Q=z^2$.  Let
\[
 F\in\{C,N\},
 \qquad b=z^kg\in\Kalg_{J-1},
 \qquad g\in R_n.
\]
The truncation conditions give $0\le k\le p-1$ and $0\le n\le2p-2$, while the contact
degree gives
\begin{equation}\label{eq:nJk}
 n=J+1-2k.
\end{equation}

\begin{lemma}[$Q\delta_t$ identity]\label{lem:QD}
For $k\ge1$ and $a\in\Kalg_1$,
\begin{equation}\label{eq:QD}
 d_3(a\wedge Q\wedge\delta_tb)
 =[a,Q]\wedge\delta_tb-[a,\delta_tb]\wedge Q
 +4k(J-k-2)b\wedge a.
\end{equation}
The first two terms lie in
$\Kalg_3\wedge\Kalg_{J-3}$ and $\Kalg_{J-2}\wedge\Kalg_2$, respectively.  Hence the
coefficient on the $\Kalg_1\wedge\Kalg_{J-1}$ component is
\[
 \gamma_{J,k}=4k(J-k-2).
\]
\end{lemma}

\begin{proof}
Since $\delta_tb=2kz^{k-1}g$, direct use of the contact bracket gives
$[Q,\delta_tb]=4k(J-k-2)b$.  Substitution in the complete
Chevalley--Eilenberg boundary gives \eqref{eq:QD}.
\end{proof}

When $\gamma_{J,k}=0$ and $n\ge2$, use the $B=zS^2V$ component.  For
$\varphi\in S^2V$ and $g_0\in R_{n-2}$,
\begin{equation}\label{eq:B-crossing}
 {[z\varphi,z^kg_0]}
 =z^{k+1}\{\varphi,g_0\}_{\spb}+(n-4)z^k\varphi g_0.
\end{equation}
On $\gamma_{J,k}=0$, the second coefficient satisfies
$n-4\equiv1-J\pmod p$.

\begin{proposition}[Bracket faces at the truncation boundary]\label{thm:higher-contraction}
Assume $J\equiv1\pmod p$. For $f\in C$, $v\in V$, and
horizontal $H$ of degree $J-2p+5$, put
\[
 T_C=f\wedge Q\wedge2(p-2)z^{p-3}H,\qquad
 T_N=zv\wedge Q\wedge2(p-2)z^{p-3}H.
\]
The bracket part of their three-boundaries is
\begin{align*}
 d_3T_C={}&8zf\wedge z^{p-3}H
 +4z^{p-3}\{f,H\}_{\spb}\wedge Q
 +12z^{p-4}fH\wedge Q-8z^{p-2}H\wedge f,\\
 d_3T_N={}&4\bigl(z^{p-2}\{v,H\}_{\spb}+z^{p-3}vH\bigr)\wedge Q
                  -8z^{p-2}H\wedge zv.
\end{align*}
For $H'$ of degree $J-2p+7$, the source
$T_N'=zv\wedge Q\wedge2(p-3)z^{p-4}H'$ has bracket part
\[
 d_3T_N'=6\bigl(z^{p-3}\{v,H'\}_{\spb}
                    +2z^{p-4}vH'\bigr)\wedge Q
                  -24z^{p-3}H'\wedge zv.
\]
With coefficients, each expression also has the three action faces
of \eqref{eq:coefficient-boundary-change}.
\end{proposition}
\begin{proof}
Substitute the contact bracket \eqref{eq:contact-bracket} in each
of the three bracket faces. The last coefficient in the general source
$a\wedge Q\wedge2kz^{k-1}H$ is $4k(J-k-2)$.
At $k=p-2,p-2,p-3$ it gives $-8,-8,-24$.
The remaining two faces give the other displayed terms.
\end{proof}

For the two lower endpoints set
$b_-=zy$, $\beta_+=zx^2$, $q_+=z^2x$.
The four bracket identities are
\begin{align*}
 d_3(f\wedge b_-\wedge\beta_+)
 ={}&(zf_x-yf)\wedge zx^2
 +(2zxf_y+x^2f)\wedge zy
 +(-2z^2x+zx^2y)\wedge f,\\
 d_3(f\wedge b_-\wedge q_+)
 ={}&(zf_x-yf)\wedge z^2x
 +(z^2f_y+2zxf)\wedge zy
 +(-z^3+z^2xy)\wedge f,\\
 d_3(zv\wedge b_-\wedge\beta_+)
 ={}&z^2v_x\wedge zx^2+(2z^2xv_y-zx^2v)\wedge zy
 +(-2z^2x+zx^2y)\wedge zv,\\
 d_3(zv\wedge b_-\wedge q_+)
 ={}&z^2v_x\wedge z^2x+(z^3v_y-z^2xv)\wedge zy
 +(-z^3+z^2xy)\wedge zv.
\end{align*}
These follow by evaluating the three brackets separately.
For $p=7$, the total-degree-ten sources can be chosen as
\[
 A_+=x^4,\quad g_+=x^3y^6,\quad H_+=x^6y^5,
 \qquad A_-=y^4,\quad g_-=x^6y^3,\quad H_-=x^5y^6.
\]
They satisfy $[A_+,g_+]=24H_+$ and $[A_-,g_-]=-24H_-$, so
\[
 d_3(a\wedge A_\pm\wedge g_\pm)
 =[a,A_\pm]\wedge g_\pm-[a,g_\pm]\wedge A_\pm
       \pm24H_\pm\wedge a.
\]
The coefficients of $a\wedge H_\pm$ are $-24,+24$, or $4,3$ modulo seven.
For $a=zv$ the other brackets are
$[zv,A]=z\{v,A\}_{\spb}+2vA$ and
$[zv,g]=z\{v,g\}_{\spb}+7vg$.

\begin{lemma}[Trivial inducing line]\label{lem:trivial-inducing}
Assume $d=0$ and $\zeta^\dagger=0$.  Then
\[
 \Kalg_jM_j(0,0)=0
 \qquad(j>0).
\]
Consequently, on the row $h=0$ (so $r=-1$), every ordinary relation of total degree
$J\ge4$ has zero induced boundary on the surviving $\Kalg_1$ quotient.
\end{lemma}

\begin{proof}
Let $u_0$ span the trivial inducing line, let $X\in\Kalg_j$, and let $m$ be a negative
PBW monomial of depth $j$.  Commute $X$ through $m$.  A term using negative degree
strictly less than $j$ leaves a positive-degree element and therefore annihilates $u_0$;
a term using exactly degree $j$ leaves an element of $\Kalg_0$, which also annihilates
$u_0$.  There is no term using more than the total depth of $m$.  Hence $Xm u_0=0$.
For a degree-$J$ relation with degree-one factor $a\in\Kalg_1$ and partner
$b\in\Kalg_{J-1}$, the only term of its coefficient boundary with head in $\Kalg_1$
is $-bm\otimes a$, with $m\in M_{J-1}(0,0)$; the preceding identity makes it zero.
\end{proof}

\subsection{The truncation obstruction}\label{sec:maxz-natural}

Let \(R=\kk[x,y]/(x^p,y^p)\), put \(\Omega=x^{p-1}y^{p-1}\),
and retain an untruncated Laurent variable \(z\).

\begin{lemma}[The principal-part obstruction]\label{lem:principal-part}
Let $p\ge5$, $0\le d<p$, and let $P$ be a Laurent tensor of
fibre $L(d)$ and any density, with pole order less than $p$.
Suppose that $L_fP$ is polynomial in $z$ for every positive contact
monomial $f$. If its reduced generator values extend to a cocycle,
then $P$ has zero principal part at $z=0$.
\end{lemma}
\begin{proof}
Let $\psi_P(f)$ be the reduction of $L_fP$ modulo $z^p$.
Reduction of polynomial tensor values commutes with the action.
Below total contact degree $2p-2$, the untruncated bracket has no
$z^p$ term, so $\psi_P$ satisfies the cocycle identity there.
By Lemma~\ref{lem:degree-one-generation}, a cocycle with the same
generator values must agree with $\psi_P$ through degree $2p-4$,
in particular at $Q=z^2$ and $z^{p-1}$.
Their truncated bracket is zero, while their untruncated bracket is
$2(p-3)z^p$. Since $L_{z^p}=2z^p\partial_z$, the defect is
\begin{equation}\label{eq:principal-part-defect}
 L_Q\psi_P(z^{p-1})-L_{z^{p-1}}\psi_P(Q)
 =-12\overline{z^p\partial_zP}
 =12\sum_{a=1}^m aA_a z^{p-a-1},
\end{equation}
where $P=P_{\rm pol}+\sum_{a=1}^m A_a z^{-a}$ and $m<p$.
The powers of $z$ are distinct and nonnegative, and each $12a$ is a unit.
The cocycle identity therefore forces every $A_a$ to be zero.
\end{proof}

For a nonzero constant fibre vector \(v\in L(d)\), consider
\(P=-\Omega v/z\) in the tensor density module of fibre character
\(\lambda=-\zeta^\dagger\).
The contact action on its positive generators gives
\begin{equation}\label{eq:top-pole-data}
 L_{C_i}P=0,\qquad
 L_{b_+}P=x^{p-1}y^{p-2}v,\qquad
 L_{b_-}P=-x^{p-2}y^{p-1}v,\qquad
 L_QP=2(3-\lambda)\Omega v .
\end{equation}
Indeed, every term in the cubic action contains a horizontal
factor multiplying \(\Omega\), or has horizontal degree above \(2p-2\).
In the natural action only \(z\partial_y\) and \(-z\partial_x\)
remain. The horizontal terms in the \(Q\)-action vanish for the
same reason.

For all positive inputs \(L_fP\) is polynomial in \(z\).
For generators of degrees one and two this follows from
\eqref{eq:top-pole-data} and their brackets. For a homogeneous input
of contact degree \(j\ge3\), a term containing \(z^{-a}\), \(a>0\),
would have horizontal degree \(j+2p-4+2a>2p-2\).
By \eqref{eq:principal-part-defect}, their defect is
\begin{equation}\label{eq:top-pole-defect}
 L_Q\psi_P(z^{p-1})-L_{z^{p-1}}\psi_P(Q)
       =-12\Omega v z^{p-2}\ne0.
\end{equation}
By Lemma~\ref{lem:principal-part}, we exclude every nonzero datum in this family.

There is a scalar version one horizontal degree lower.
For \(H\in R_m\), \(m=2p-3\) or \(2p-2\), set
\(\lambda=1-m\). In the Laurent density module,
\begin{equation}\label{eq:scalar-pole-data}
 L_{z^kF_n}(H/z)
   =z^{k-1}\{F_n,H\}_{\rm sp}
       +(k+n-2)z^{k-2}F_nH,\qquad F_n\in R_n .
\end{equation}
These values are polynomial for every positive input.
For \(k=1,n=1\) the second coefficient vanishes; for \(k=1,n\ge2\)
horizontal truncation kills its numerator. At \(k=0,n\ge3\)
the products vanish, and so do the Poisson terms:
the only endpoint to check is
\(\{R_3,R_{2p-3}\}_{\rm sp}=0\), obtained by differentiating
\(x^{p-1}y^{p-2}\) and \(x^{p-2}y^{p-1}\).
For example the possible top coefficient in
\(\{x^2y,x^{p-2}y^{p-1}\}_{\rm sp}\) is \(p=0\).
The reduced cochain has values
\[
 \chi(C_i)=0,\qquad \chi(b_+)=H_y,\qquad
 \chi(b_-)=-H_x,\qquad \chi(Q)=0,
\]
and defect \(12Hz^{p-2}\) by \eqref{eq:principal-part-defect}.
Applying Lemma~\ref{lem:principal-part}, we obtain \(H=0\) for a cocycle with these data.

\section{Restricted cocycles and coefficient boundaries}\label{sec:restricted}

For a left restricted \(\mathfrak n\)-module \(T\), an ordinary
one-cocycle \(\chi\) is restricted precisely when
\begin{equation}\label{eq:restricted-cocycle-test}
 \chi(a^{[p]})=a^{p-1}\chi(a)\qquad(a\in\mathfrak n).
\end{equation}
This follows by pairing with the \(p\)-power differential
\eqref{eq:res-p}. Inner cocycles satisfy
\eqref{eq:restricted-cocycle-test} by the module identity
\(\rho(a^{[p]})=\rho(a)^p\).
Thus ordinary relations give upper bounds for the restricted groups;
to attain such a bound it suffices to construct the indicated
restricted cocycles.

We first check \eqref{eq:restricted-cocycle-test} for connecting cocycles.
If
\(0\to T\to\widetilde T\to Q_0\to0\)
is a sequence of restricted modules and \(q\in Q_0^{\mathfrak n}\),
the connecting cocycle \(a\mapsto a\widetilde q\) is restricted,
since
\[
 \chi(a^{[p]})=a^{[p]}\widetilde q
             =a^p\widetilde q=a^{p-1}\chi(a).
\]
For closed differential forms on the truncated function algebra,
the identity \eqref{eq:Cartier-restricted-identity} below gives the
same condition directly.
The tensor and divided-power constructions in
Section~\ref{sec:first-two} and the closed-form construction below
supply these lower bounds.

Let $A=\kk[x,y,z]/(x^p,y^p,z^p)$, with contact action
\[
 D_f=(-f_y+xf_z)\partial_x+(f_x+yf_z)\partial_y
                +(2f-xf_x-yf_y)\partial_z.
\]
Put $\theta=dz-y\,dx+x\,dy$ and $H^*=\Omega_A^1/A\theta$.
Let $\bar d:A\to H^*$ be exterior differentiation followed by projection;
in the frame $X=dx$, $Y=dy$ it is
\[
 \bar d g=(g_x+yg_z)X+(g_y-xg_z)Y.
\]
We use the function module $A=T(0,0)$ and the horizontal form module
$H^*=T(1,1)$, as defined in Section~\ref{sec:first-two}.

\begin{lemma}[The contact multiplier]\label{lem:contact-multiplier}
For every $p\ge5$, the maps
\[
 c(f)=f_z,\qquad c_H(f)=\bar d(f_z),\qquad f\in\mathfrak n,
\]
are restricted cocycles with values in $A$ and $H^*$, respectively.
They give nonzero horizontal invariant classes of internal degrees $0$ and $-1$.
\end{lemma}
\begin{proof}
The contact Lie derivative satisfies $\mathcal L_{D_f}\theta=2f_z\theta$.
Its commutator identity on $A\theta$ shows that $c$ is an ordinary cocycle.
On this line the operator is $D_f+M_{2f_z}$, with $M_a$ multiplication by $a$.
Using $\mathcal L_{D_f}^p=\mathcal L_{D_f^p}$ and
\[
 (D+M_a)^p=D^p+M_{D^{p-1}a+a^p},
\]
we obtain $c(f^{[p]})=D_f^{p-1}c(f)+c(f)^p$.
For $f\in\mathfrak n$, the function $f_z$ lies in the augmentation ideal
of $A$, so $c(f)^p=0$. Thus $c$ is restricted.
The map $\bar d$ is equivariant, so $c_H$ is restricted as well.
Their generator values are
\[
 c(C_i)=c_H(C_i)=0,\quad
 (c(b_+),c(b_-))=(x,y),\quad
 (c_H(b_+),c_H(b_-))=(X,Y).
\]
On an input of contact degree $j$, their coefficient degrees are $j$ and $j-1$.
A homogeneous primitive for $c$ would be constant and have zero coboundary;
a primitive for $c_H$ would have negative coefficient degree.
Neither cocycle is inner. The horizontal group preserves $\theta$ and
commutes with $\partial_z$ and $\bar d$, so both classes are invariant.
\end{proof}

\subsection{Cartier cocycles}\label{sec:terminal-proof}

\begin{proof}[Proof of Theorem~\ref{thm:qtop}]
There is a \(K(3;\mathbf1)^{(1)}\)-module isomorphism
\(M(0,0)\cong A^*\) sending the inducing vector to evaluation at zero.
Indeed, the negative fields \(D_x,D_y,D_1\) are independent at the origin.
In increasing differential order the PBW distributions have leading
symbols with nonzero coefficients \(a!b!c!2^c\), \(0\le a,b,c<p\).
They therefore give a triangular basis of \(A^*\).
We use the resulting pairing
\[
 \langle u_0,F\rangle=F(0),\qquad
 \langle am,F\rangle=-\langle m,D_aF\rangle .
\]

Contraction of the closed forms \(x^{p-1}dx,y^{p-1}dy\) gives
\[
 \chi_x(f)=-x^{p-1}f_y,\qquad
 \chi_y(f)=y^{p-1}f_x.
\]
Closedness is the ordinary cocycle identity.
For a derivation \(D\) of a commutative algebra of characteristic \(p\),
\begin{equation}\label{eq:Cartier-restricted-identity}
 D^{p-1}(u^{p-1}Du)=u^{p-1}D^pu-(Du)^p.
\end{equation}
To verify this, use the multiplication operator \(M_b\) and expand
\((D+M_b)^p=D^p+M_{D^{p-1}b+b^p}\).
The mixed Bell coefficients are divisible by \(p\).
Set \(g=1+\nu u\), \(b=Dg/g\), use
\(D+M_b=M_g^{-1}DM_g\), and compare the coefficient of \(\nu^p\).
For \(f\in\mathfrak n\), the values \(D_fx,D_fy\) vanish at the origin,
so their \(p\)-th powers are zero in \(A\).
Equation~\eqref{eq:Cartier-restricted-identity} proves
\[
 D_f^{p-1}\chi_i(f)=\chi_i(f^{[p]}),\qquad i=x,y.
\]
Thus both cocycles are restricted.

If \(a\chi_x+b\chi_y\) is inner, its primitive \(F\) may be taken of
weighted degree \(p\).
The equation \(D_{z^2}F=0\), with
\[
 D_{z^2}=2z(x\partial_x+y\partial_y)+2z^2\partial_z,
\]
forces \(F\) to be horizontal: on
\(z^kF_{p-2k}(x,y)\) it has coefficient \(-2k\).
For \(1\le k\le(p-1)/2\) these coefficients are units and
\(z^{k+1}\ne0\).
The equations for \(zx,zy\) then give
\(F_x=ax^{p-1}\), \(F_y=by^{p-1}\).
Neither derivative can contain its displayed pure top monomial in the
truncated horizontal ring. Hence \(a=b=0\).

For horizontal functions the pairing is
\[
 \langle x^ay^bu_0,F\rangle
       =(-1)^a(\partial_y^a\partial_x^bF)(0).
\]
The two cubic cycles in the theorem consequently have pairing matrix
\(\operatorname{diag}(6,-6)\) with \((\chi_x,\chi_y)\).
They are cycles because their coefficient action contains \(x^p\) or \(y^p\).
The two truncation cycles are cycles by
\([b_+,t]=-2x\), \([b_-,t]=-2y\), and their nonzero pairings are
\[
 \langle x^{p-1}tu_0,\chi_y(b_+)\rangle=2,\qquad
 \langle y^{p-1}tu_0,\chi_x(b_-)\rangle=-2 .
\]
Their internal degree is \(p\), so \(h=p+1\);
PBW duality gives \(q=3p-4\), \(\zeta=p-4\).

The integer horizontal weights of the closed forms modulo exact forms
are \(p,-p\).
Expanding \(x\mapsto x+\nu y\) gives the root-group coefficient \(\nu^p\)
on their classes; the other root group is symmetric.
They form \(V^{(1)}\), and the injective contraction map preserves this action.
\end{proof}

\subsection{The three truncated-Koszul classes}\label{sec:truncated-koszul}

For each of the three $t$-positive classes in Proposition~\ref{prop:Koszul}, the leading
degree-three relation admits a canonical lift to the kernel of the Lie boundary.

Put
\[
 \beta_0=zxy,
 \qquad
 \beta_+=zx^2,
 \qquad
 \gamma_0=x^3,
 \qquad
 \gamma_1=x^2y,
 \qquad
 \alpha_0=x^2y^2.
\]
The leading $z$-degree-two relation of weight $+1$ is
\[
 \rho_+^{(2)}=2b_+\wedge \beta_0+b_-\wedge \beta_+.
\]
Its Lie boundary is
\[
 d_{\mathrm{Lie}}\rho_+^{(2)}=3zx^2y.
\]
Consider the correction subspace
\[
 \mathscr E_+=\kk(b_+\wedge \alpha_0)\oplus\kk(\gamma_1\wedge \beta_0).
\]
With respect to the basis $(zx^2y,x^3y^2)$, the bracket map on $\mathscr E_+$ is
\begin{equation}\label{eq:obstruction-matrix}
 \begin{pmatrix}2&1\\2&-1\end{pmatrix},
 \qquad \det=-4.
\end{equation}
Since this matrix is invertible, there is a unique element of $\mathscr E_+$ whose sum with
$\rho_+^{(2)}$ lies in $\ker d_{\mathrm{Lie}}$.

\begin{proposition}[Canonical lift of $\rho_+^{(2)}$]\label{prop:Rplus}
The element
\begin{equation}\label{eq:Rplus}
 R_+
 =2b_+\wedge \beta_0+b_-\wedge \beta_+
 -\frac34b_+\wedge \alpha_0
 -\frac32\gamma_1\wedge \beta_0
\end{equation}
belongs to $\ker d_{\mathrm{Lie}}$.  Its coefficients are uniquely determined by
\eqref{eq:obstruction-matrix}.
\end{proposition}

\begin{proof}
The required brackets are
\[
 {[b_+,\beta_0]}=z^2x+zx^2y,
 \qquad
 {[b_-,\beta_+]}=-2z^2x+zx^2y,
\]
\[
 {[b_+,\alpha_0]}=2zx^2y+2x^3y^2,
 \qquad
 {[\gamma_1,\beta_0]}=zx^2y-x^3y^2.
\]
The sum of the first two contributions has boundary $3zx^2y$.  For the coefficients
$\nu_1$ of $b_+\wedge\alpha_0$ and $\nu_2$ of $\gamma_1\wedge\beta_0$, the equations
\[
 3+2\nu_1+\nu_2=0,
 \qquad
 2\nu_1-\nu_2=0
\]
have the unique solution $\nu_1=-3/4$ and $\nu_2=-3/2$, yielding \eqref{eq:Rplus}.
\end{proof}

Take a coefficient on the highest-weight line,
\[
 w=x^{p-1}t^{c+1}\otimes u_{\mathrm{hw}},
 \qquad W=\beta_0w.
\]
On this line
\[
 b_+w=\gamma_1w=\alpha_0w=\beta_+w=0.
\]
Using the left homology convention
\[
 d_2(m\otimes a\wedge b)=am\otimes b-bm\otimes a+m\otimes[a,b],
\]
one obtains
\begin{equation}\label{eq:Rplus-boundary}
 d_2(w\otimes R_+)
 =-2W\otimes b_+ +(b_-w)\otimes \beta_+ +\frac32W\otimes \gamma_1.
\end{equation}
The $\beta_+$ component is eliminated by the bracket section
\begin{equation}\label{eq:Bplus-section}
 U=\frac14(\gamma_0\wedge b_- - \gamma_1\wedge b_+),
 \qquad
 {[\gamma_0,b_-]}-[\gamma_1,b_+]=4\beta_+.
\end{equation}
If $m=b_-w$, then
\[
 {[\gamma_1,b_-]}=2\beta_0-\alpha_0,
 \qquad \gamma_1m=2W.
\]
After the substitution, the $b_+$ coordinate acquires the additional term
$\frac12W$.

\begin{theorem}[Reduction of the $R_+$ boundary]\label{thm:frob-scalar}
For $d=0$ and $(p,c)=(5,1),(7,1),(7,2)$, the boundary of
$x^{p-1}t^{c+1}u_0\otimes R_+$, after bracket-section elimination,
belongs to the image of the quadratic relation boundaries.
Its induced map on this quadratic quotient is zero.
\end{theorem}
\begin{proof}
First fix the raw coordinate for a general inducing weight $d$.
For $u_{\rm hw}=X^d$ and $\eta_+=x^{p-1}t^cu_{\rm hw}$,
commuting $\beta_0=zxy$ through $x^at^n$ gives
\[
 \beta_0x^at^nu_{\rm hw}=2n(d+a)x^at^{n-1}u_{\rm hw}.
\]
Here $[\beta_0,t]=-2xy$, $[\beta_0,x]=-b_+-\gamma_1$, and
the $\gamma_1$ terms vanish on this highest-weight line.
The term from $\beta_0t^n$ contributes $2nd$; the $a$ commutator
terms contribute $2n$ each. Hence, for $w=x^{p-1}t^{c+1}u_{\rm hw}$,
\begin{equation}\label{eq:lambda-c}
 W=\beta_0w=2(c+1)(d-1)\eta_+.
\end{equation}

Now set $d=0$. Eliminating the $\beta_+$ head in
\eqref{eq:Rplus-boundary} by \eqref{eq:Bplus-section} gives
\[
 R_{\rm red}=-\tfrac32W\otimes b_+
       +\tfrac14 b_-^2w\otimes C_0
       +(\tfrac32W-\tfrac14b_+b_-w)\otimes C_1
       -\tfrac14 C_0b_-w\otimes b_-.
\]
The last term is zero on this weight-$p$ line.
Put
\begin{align*}
 A&=C_0\otimes x^{p-1}y^2t^{c-1}u_0,&
 H&=b_+\otimes x^{p-1}t^cu_0,\\
 B&=C_0\otimes x^{p-2}yt^cu_0,&
 C&=C_1\otimes x^{p-1}t^cu_0,\\
 D&=C_0\otimes x^{p-3}t^{c+1}u_0.&&
\end{align*}
These are all the possible weight-$p$ coordinates at this depth.
Indeed, the heads have weights $3,1,-1,-3$; the last two require
horizontal coefficient weight exceeding the truncation range.
For head weight three, the three possible $y$-exponents are $2,1,0$;
the next one would require $x$-exponent at least $p$.

Twice using $[b_-,t]=-2y$ gives the lowest component
\[
 b_-^2w=4c(c+1)x^{p-1}y^2t^{c-1}u_0+\text{terms of $t$-exponent at least $c$}.
\]
Thus the $(A,H)$ coordinates of $R_{\rm red}$ are
$(c(c+1),3(c+1))$.
For $r_4=C_0\wedge C_1+3C_0\wedge b_+$, let
\[
 q_0=\partial_2(x^{p-2}y^2t^cu_0\otimes r_4).
\]
Its $(A,H)$ coordinates are $(6c,18)$: the first is the one-step
$t$-drop of $-3b_+$, and the second is the coefficient in
$3C_0(x^{p-2}y^2t^cu_0)=18x^{p-1}t^cu_0$.
Subtracting $(c+1)q_0/6$ therefore removes both $A$ and $H$.

Let $f=-\operatorname{ad}(y^2)/2$ and take the three remaining quadratic sources
\[
 x^{p-3}yt^{c+1}u_0\otimes r_4,\qquad
 x^{p-2}t^{c+1}u_0\otimes fr_4,\qquad
 x^{p-4}t^{c+2}u_0\otimes r_4,
\]
where $fr_4=2C_0\wedge C_2+3C_0\wedge b_-+9C_1\wedge b_+$.
Their boundaries have zero $A,H$ coordinates, and on $(B,C,D)$ they satisfy
\[
 q_1\equiv6(c+1)B,\qquad
 q_2\equiv6(c+1)B+18(c+1)C\pmod{\kk D},\qquad
 q_3=6(c+2)D.
\]
The third boundary spans $\kk D$.
The coefficients follow from the single $t$-drop in the natural heads;
higher $t$-exponents have insufficient horizontal weight.
All three pivots are units at the stated pairs, and all source
$t$-exponents are less than $p$. They remove the remaining $B,C,D$
coordinates, proving the claim.
\end{proof}

\section{Coupled horizontal relation coordinates}\label{sec:tfree-horizontal}

\begin{lemma}[Natural projection]\label{lem:c0-central}
For $a\in R_h\otimes L(d)$ and $b\in N$, the natural $t$-free
coordinate of a coefficient boundary is
\begin{equation}\label{eq:c0-central}
 \operatorname{pr}_{N,c=0}\partial_2(ta\otimes b\wedge Q)
       =4(\zeta^\dagger-h)a\otimes b.
\end{equation}
\end{lemma}
\begin{proof}
Use $[b,Q]=0$ and $Q(ta)=t(Qa)-4za$.
The central element $z$ acts on the horizontal coefficient by
$\zeta^\dagger-h$. This gives the displayed natural coordinate.
For example, at $h=0$ the full boundary is
\[
 \partial_2(tu\otimes b_+\wedge Q)
       =-2xu\otimes Q+4\zeta^\dagger u\otimes b_+.
\]
\end{proof}

For a restricted $\slie$-module $W$, define
\begin{equation}\label{eq:PhiW}
 \Phi_W:S^3V\otimes W\longrightarrow V\otimes W,\qquad
 \Phi_W(f\otimes w)=x\otimes\rho_W(f_x)w+y\otimes\rho_W(f_y)w.
\end{equation}

\begin{proposition}[The projected natural operator]\label{prop:c0-natural}
For $0\le s,d<p$,
\begin{equation}\label{eq:Phi-coker}
 \coker\Phi_{L(s)\otimes L(d)}\cong
 \begin{cases}
 L(1),&s=d,\\
 L(0),&|s-d|=1,\\
 0,&|s-d|\ge2.
 \end{cases}
\end{equation}
\end{proposition}
\begin{proof}
Write $E=X\partial_Y$, $F=Y\partial_X$, $H=X\partial_X-Y\partial_Y$.
The four monomial images are
\[
 6x\otimes Ew,\quad -2x\otimes Hw+2y\otimes Ew,\quad
 -2x\otimes Fw-2y\otimes Hw,\quad -6y\otimes Fw.
\]
The transpose kernel is therefore the space of $(a,b)\in W^*\oplus W^*$ with
\[
 Ea=Fb=0,\qquad Eb=Ha,\qquad Fa=-Hb.
\]
In a toral weight component put $Ha=\alpha a$, $Hb=(\alpha-2)b$.
Then $\alpha a=EFa=\alpha(2-\alpha)a$, hence
$\alpha(1-\alpha)a=0$.
For $\alpha=0$ one has $b=0$ and $a$ invariant;
for $\alpha=1$ the pair is the image of a homomorphism $V\to W^*$.
The solutions with $a=0$ give the second invariant coordinate.
Thus
\[
 \ker\Phi_W^*\cong
 V^*\otimes(W^*)^{\slie}\oplus\Hom_{\slie}(V,W^*).
\]
For $W=L(s)\otimes L(d)$, self-duality gives the first multiplicity
$\delta_{s,d}$. If $s\le p-2$, multiplication and its gradient section
split $V\otimes L(s)$ as $L(s+1)\oplus L(s-1)$, since $s+1$ is a unit.
Tensor duality treats the case $s=p-1,d<p-1$ by interchanging $s,d$.
For $s=d=p-1$, the multiplication sequence
\[
 0\longrightarrow L(p-2)\longrightarrow V\otimes L(p-1)
     \longrightarrow S^pV\longrightarrow0
\]
shows that the square of the Casimir
$\Omega=H^2+2H+4FE$ kills its middle term:
the Casimir is zero on both endpoints, whereas it is $-1$ on the
proposed target $L(p-1)$.
It follows that the second multiplicity is one for $|s-d|=1$ and zero
otherwise. The one-dimensional multiplicity spaces have trivial rational
$SL_2$-action. Dualizing proves the formula.
\end{proof}

Set $A_0=x^4$, $\beta_+=zx^2$ and
$\rho_5=2C_0\wedge\beta_++b_+\wedge A_0$.
Using the bracket sections
\[
 \beta_+=([C_0,b_-]-[C_1,b_+])/4,\qquad A_0=-[C_0,b_+],
\]
the negative of its reduced boundary has natural and cubic coordinates
\begin{equation}\label{eq:full-horizontal-graph}
\begin{aligned}
 G_N(m)={}&b_+\otimes(b_+C_0-2C_0b_+-C_1C_0/2)m
                   +b_-\otimes C_0^2m/2,\\
 G_C(m)={}&C_0\otimes(b_+^2+2\beta_+-b_-C_0/2)m
                   +C_1\otimes b_+C_0m/2.
\end{aligned}
\end{equation}
Indeed, substitute
$\alpha_{[a,b]}(m)=bm\otimes a-am\otimes b$
for each bracket head in the coefficient boundary.
The $+2\beta_+$ term is one of the resulting cubic coordinates.

There is also a return across the positive $t$-layers in $G_N$.
On symmetric PBW coordinates let
\[
 \mathscr S=\sum_{j=0}^{p-1}\frac{(-t/2)^j}{j!}
                       (\partial_x\partial_y)^j.
\]
For horizontal inputs write
\begin{align*}
 A&=3x\partial_y^2+6E\partial_y-\tfrac t2\partial_y^3=A_0'+tA_1',\\
 B&=y\partial_y^2-2x\partial_x\partial_y-2H\partial_y-2E\partial_x
                       +\tfrac t2\partial_x\partial_y^2,\\
 T&=-2x\partial_t+
 (\zeta^\dagger+H-x\partial_x-y\partial_y-t\partial_t)\partial_y-2E\partial_x.
\end{align*}
These are $\mathscr S^{-1}L_{C_0}\mathscr S$,
$\mathscr S^{-1}L_{C_1}\mathscr S$ and
$\mathscr S^{-1}L_{b_+}\mathscr S$, respectively.
Their formulas follow by commuting with the negative left multiplications
$x+(t/2)\partial_y$, $y-(t/2)\partial_x$ and $t$, and evaluating on the inducing fibre.
If $T_h,B_h$ denote their horizontal parts, the two natural coordinates are
\[
 K_x=T_hA_0'-2A_0'T_h-B_hA_0'/2+x\partial_y^3,
 \qquad K_y=(A_0')^2/2.
\]
The last term in $K_x$ is the product of $-2x\partial_t$ in $T$
with $-t\partial_y^3/2$ in $A$.

\subsection{Normalized natural readings}\label{sec:normalized-natural}

We fix the natural-coordinate readings of
\eqref{eq:full-horizontal-graph}.
Put \(s=\min(h,2p-2-h)\), and use
\[
 \iota_h(x^{s-k}y^k)=
 \begin{cases}
 x^{s-k}y^k,&h<p,\\
 (-1)^k\binom sk^{-1}x^{p-1-k}y^{h-p+1+k},&h\ge p .
 \end{cases}
\]
The apolar pairing is
\[
 \beta_d(x^{d-k}y^k,V_j)
   =\delta_{k+j,d}\frac{(-1)^k}{\binom dk},
 \qquad V_j=X^{d-j}Y^j .
\]
The \(E,F\) recursions and the binomial identity show that this
pairing is invariant. Its normalization is
\[\beta_d(x^d,Y^d)=1 .\]
For \(u=\iota_h^{-1}(f)\), define the natural readings
\[
\begin{array}{c|cc}
 &b_+\otimes f\otimes v&b_-\otimes f\otimes v\\ \hline
\ell_-\ (d=s-1)&\beta_{s-1}(\partial_yu,v)&-\beta_{s-1}(\partial_xu,v)\\
\ell_0\ (d=s)&\beta_s(u,v)&0\\
\ell_+\ (d=s+1)&\beta_{s+1}(xu,v)&\beta_{s+1}(yu,v).
\end{array}
\]
Here \(\ell_0\) reads the \(b_+\)-coordinate of the invariant
natural quotient. Substitution in the four monomial images of
\(\Phi\) shows that these functionals annihilate that image.
They take value one on the respective representatives
\[
\begin{aligned}
 r_-&=b_+\otimes\iota_h(x^{s-1}y)\otimes Y^{s-1},\\
 r_0&=b_+\otimes\iota_h(x^s)\otimes Y^s,\\
 r_+&=b_+\otimes\iota_h(x^s)\otimes Y^{s+1}.
\end{aligned}
\]

Take a symmetric-coordinate monomial \(\varphi=x^ay^bV_i\),
with \(a+b=h+2\), whose PBW representative is \(\mathscr S\varphi\).
The \(K_x\)-monomials and \(K_y\)-monomials are respectively
\[
\begin{aligned}
 &(x^{a+1}y^{b-3}V_i,\ x^ay^{b-2}V_{i-1},\
                                  x^{a-1}y^{b-1}V_{i-2}),\\
 &(x^{a+2}y^{b-4}V_i,\ x^{a+1}y^{b-3}V_{i-1},\
                                  x^ay^{b-2}V_{i-2}).
\end{aligned}
\]
Their coefficients, with out-of-range monomials zero, are
\begin{equation}\label{eq:natural-six-coefficients}
\begin{aligned}
 t_0&=\tfrac12b(b-1)(b-2)(12a+3b-6\zeta^\dagger+17),\\
 t_1&=3bi(b-1)(7a+b-2\zeta^\dagger+9),&
 t_2&=18abi(i-1),\\
 u_0&=\tfrac92b(b-1)(b-2)(b-3),\\
 u_1&=18b(b-1)(b-2)i,&u_2&=18b(b-1)i(i-1).
\end{aligned}
\end{equation}
For an explicit expansion, put
\(u_*=\zeta^\dagger+d-2i-a-b+1\).
The factor multiplying \(b(b-1)(b-2)\) in \(t_0\) is
\[
 -3u_*+3a+3(d-2i)-3b/2+23/2 .
\]
The factor multiplying \(3bi(b-1)\) in \(t_1\) is
\(-2u_*+5a-b+2(d-2i)+11\).
Substitution gives the two displayed formulas.
The coefficients \(u_j\) follow from the four terms of
\((A'_0)^2/2\); adding its two mixed terms, which include the
\(x\partial_y^3\) return in \(K_x\), gives the factor \(18\).

On the high side \(h=2p-2-s\), \(4\le s\le p-1\),
use the sources
\[
\begin{aligned}
 m_-&=-\frac{(s-1)!}{s}
       \mathscr S(x^{p-3}y^{p-s+3}Y^{s-1}),\\
 m_0&=s!\mathscr S(x^{p-2}y^{p-s+2}Y^s),\\
 m_+&=(s+1)s!\mathscr S(x^{p-2}y^{p-s+2}Y^{s+1}),
\end{aligned}
\]
for the legal fibre degrees \(d=s-1,s,s+1<p\).
For the unscaled minus source, the six apolar weights are all
\(-s\). For the zero source, the three \(K_x\)-weights are one
and the \(K_y\)-weights are zero. For the plus source the
\(K_x\)-weights are
\(1,s/(s+1),(s-1)/(s+1)\), while the last two \(K_y\)-weights
are \(-1/(s+1),-2/(s+1)\); its first \(K_y\)-term is truncated.
Inserting \(\zeta^\dagger=h=-2-s\) in
\eqref{eq:natural-six-coefficients} gives the unscaled readings
\[
\begin{aligned}
 -2s(s-3)(s-2)(s-1)(3s+10),\\
 \tfrac12s(s-1)(s-2)(3s+67),\\
 \tfrac12s(s-1)(s-2)(3s+31).
\end{aligned}
\]
For the first, divide out \((s-3)(s-2)(s-1)\):
the sums of the \(K_x,K_y\) terms are \(3s/2+38,9s/2-18\),
then multiply by \(-s\).
For the second, the remaining sum is
\(-(3s+11)/2+3(s+1)+36\).
For the third, the last two corresponding fibre terms cancel,
leaving \(-(3s+11)/2+3(s+1)+18\).
Multiplication by the three source factors therefore gives
\begin{equation}\label{eq:corrected-natural-kappas}
\begin{aligned}
 \widehat\kappa_-&=2(s-3)(s-2)^2(s-1)^2(3s+10)(s-3)!,\\
 \widehat\kappa_0&=\tfrac12s^2(s-1)^2(s-2)(3s+67)(s-2)!,\\
 \widehat\kappa_+&=\tfrac12s^2(s-1)^2(s-2)(s+1)(3s+31)(s-2)! .
\end{aligned}
\end{equation}
The common unit is one for the sign in
\eqref{eq:full-horizontal-graph}.
For instance, the legal high-side evaluations are
\[
\begin{array}{cc|ccc}
p&s&\widehat\kappa_-&\widehat\kappa_0&\widehat\kappa_+\\ \hline
5&4&4&2&\text{--}\\
7&4&2&2&5\\
7&5&2&3&6\\
7&6&0&3&\text{--}.
\end{array}
\]

The four horizontal heads of an ordered
one-\(t\) source \(tx^ay^bY^d\) can be read directly from the same graph.
Let \(D_f\) be the ordered horizontal action.
The commutators with \(t\) give
\[
\begin{aligned}
 \pi_0G_N(t\varphi)&=b_+\otimes2xD_{C_0}\varphi,\\
 \pi_0G_C(t\varphi)&=
 C_0\otimes(-4xD_{b_+}-2D_{x^2}+yD_{C_0})\varphi
                  -C_1\otimes xD_{C_0}\varphi .
\end{aligned}
\]
Here
\[
 D_{C_0}=3x\partial_y^2+6E\partial_y,\quad
 D_{x^2}=2x\partial_y+2E,\quad
 D_{b_+}=(\zeta^\dagger+H-2x\partial_x-y\partial_y)\partial_y
                                                     -2E\partial_x .
\]
The last formula follows by transporting \(T\) through
\(\mathscr S^{-1}\varphi=\varphi+t\partial_x\partial_y\varphi/2+\cdots\).
Its additional contribution is \(-x\partial_x\partial_y\).

Write the four heads as
\[
\begin{aligned}
 \pi_0G_N&=b_+\otimes
 (n_1x^{a+1}y^{b-1}V_{d-1}+n_0x^{a+2}y^{b-2}V_d),\\
 \pi_0G_C&=C_0\otimes
 (c_0x^{a+1}y^{b-1}V_d+c_1x^ay^bV_{d-1})\\
 &\quad+C_1\otimes
 (e_1x^{a+1}y^{b-1}V_{d-1}+e_0x^{a+2}y^{b-2}V_d).
\end{aligned}
\]
Substitution gives
\begin{equation}\label{eq:one-t-six-coefficients}
\begin{aligned}
 n_1&=12bd,&n_0&=6b(b-1),\\
 c_0&=b(8a+7b+4d-4\zeta^\dagger-11),&
 c_1&=d(8a+6b-4),\\
 e_1&=-6bd,&e_0&=-3b(b-1).
\end{aligned}
\end{equation}
In particular,
\[
 D_{b_+}\varphi
 =b(\zeta^\dagger-d-2a-b+1)x^ay^{b-1}V_d
                -2ad\,x^{a-1}y^bV_{d-1},
\]
which gives \(c_0,c_1\) upon insertion in the \(C_0\)-head.
At the central density \(\zeta^\dagger=h\) the evaluations are
\[
\begin{array}{ccc|cc|cc|cc|cc|c}
p&h&d&a&b&n_1&n_0&c_0&c_1&e_1&e_0&\ell\\ \hline
5&4&3&1&3&3&1&2&1&1&2&4\\
7&7&5&4&3&5&1&4&6&1&3&6\\
7&9&4&4&5&2&1&5&1&6&3&6\\
7&6&5&3&3&5&1&6&1&1&3&6\\
7&5&4&2&3&4&1&3&1&5&3&6 .
\end{array}
\]
The respective apolar weights on \(n_0,n_1\) are
\((1,1),(1,1),(1,6),(1,1),(1,3)\).
Thus the last column is
\(1+3,1+5,1+12,1+5,1+12\), in the indicated characteristics.
For the equivariant sources \(t\mathscr S(x^ay^b)Y^d\),
replace \(c_0\) by \(c_0-4ab\).
Indeed, the next source term
\(-ab\,t^2x^{a-1}y^{b-1}Y^d/2\) returns through the two
\(t\)-drops in \(b_+^2\), with coefficient \(8x^2\).
The five corrected \(c_0\)'s are \(0,5,2,5,0\).

Finally, for \(p=7\), the low-side horizontal readings obtained from
\eqref{eq:natural-six-coefficients}, with the same apolar functionals,
are as follows:
\[
\begin{array}{cc|cc|ccc|ccc|cc|c|c}
h&d&a&b&t_0&t_1&t_2&u_0&u_1&u_2&\ell_x&\ell_y&
 \text{source factor}&\ell\\ \hline
4&3&2&4&5&1&3&3&1&1&5&0&2&3\\
4&4&3&3&2&1&5&\text{--}&\text{--}&\text{--}&2&0&3&6\\
4&5&3&3&2&3&6&0&1&4&2&3&1&5\\
5&4&3&4&0&5&2&3&6&2&2&4&5&2\\
5&5&4&3&6&5&1&\text{--}&\text{--}&\text{--}&3&0&1&3\\
5&6&4&3&6&6&5&0&4&6&3&3&6&1 .
\end{array}
\]
Here the sources are \(\mathscr S(x^ay^bY^d)\), multiplied by
the factors \(-(s-1)!/s,s!,(s+1)s!\) in their three branches,
with \(s=h\).
The last column is the source factor times \(\ell_x+\ell_y\);
for example its last entry is \(6(3+3)=1\).
The complete cohomology calculations below impose both the
natural and cubic heads of \eqref{eq:full-horizontal-graph}
and the remaining relations.

\subsection{Low-degree ordinary boundaries}\label{sec:lowrows}

We fix the negative PBW order \(x,y,t\), with \(yx=xy-t\), and
write \(u_i=X^{d-i}Y^i\) for the inducing fibre.
For \(e=\operatorname{ad}(x^2)/2\), \(f=-\operatorname{ad}(y^2)/2\), their coefficient actions are
\begin{align*}
 f(x^ay^bt^cu_i)
 &=a x^{a-1}y^{b+1}t^cu_i
   -\tfrac12a(a-1)x^{a-2}y^bt^{c+1}u_i
   +(d-i)x^ay^bt^cu_{i+1},\\
 e(x^ay^bt^cu_i)
 &=b x^{a+1}y^{b-1}t^cu_i
   -\tfrac12b(b-1)x^ay^{b-2}t^{c+1}u_i
   +i x^ay^bt^cu_{i-1}.
\end{align*}
Terms outside the exponent or fibre bounds are zero.
Put
\[
 r_4=C_0\wedge C_1+3C_0\wedge b_+,\qquad
 \rho_3=C_0\wedge Q+2b_+\wedge\beta_+.
\]
Their brackets vanish by \eqref{eq:contact-bracket}.
On both relation and coefficient factors we use \(f^{(j)}=f^j/j!\);
all indices needed below are at most four.
If \(r_a=f^{(a)}r_4\) and \(w_b=f^{(b)}w_0\), equivariance gives
\begin{equation}\label{eq:boundary-divided-recursion}
 \partial_2(w_b\otimes r_a)=
 \sum_{j=0}^a(-1)^j\binom{b+j}{j}
        f^{(a-j)}\partial_2(w_{b+j}\otimes r_4).
\end{equation}
The same formula holds with \(\rho_3\). Coefficient vectors beyond their
displayed highest-weight chains are zero.

\begin{proposition}\label{prop:low-horizontal}
Under \(\zeta^\dagger=h\), the ordinary relation quotient removes the
following components:
\[
 (h,d)=(0,3):L(0),\quad (1,1):2L(1),\quad
 (2,1):2L(0),\quad (3,0):L(0).
\]
The \((1,2)\) row is zero as well.
The \((2,1)\) entry refers to its two-dimensional scalar cycle kernel.
\end{proposition}
\begin{proof}
At \((h,d)=(0,3)\), take \(m=yY^3\).
The three actions are
\[
 C_0m=18XY^2,\qquad C_1m=6Y^3,\qquad
 b_+m=(\zeta^\dagger-3)Y^3 .
\]
Therefore
\[
 \partial_2(m\otimes r_4)=18XY^2\otimes C_1+
       3(1-\zeta^\dagger)Y^3\otimes C_0+54XY^2\otimes b_+ .
\]
The cubic apolar functional
\(\tau(C_i\otimes X^iY^{3-i})=(-1)^i/\binom3i\)
takes value \(-3(\zeta^\dagger+1)\), a unit at \(\zeta^\dagger=0\).
The other Weyl factors of \((C\oplus N)\otimes L(3)\) have highest
weights \(2,4,6\), whose Casimir values are nonzero at both primes.
Thus this boundary removes the unique scalar block.

For \((h,d,\zeta^\dagger)=(1,1,1)\), use the ordered one-chain vectors
\[
 (\mathbf A,\mathbf B,\mathbf C,\mathbf D,\mathbf E,\mathbf F,\mathbf G)
 =(C_0\otimes yY,C_1\otimes xY,C_1\otimes yX,C_2\otimes xX,
   b_+\otimes xY,b_+\otimes yX,b_-\otimes xX).
\]
The weight-one highest vectors are
\[
 c=\mathbf A-\mathbf B-\mathbf C+\mathbf D,\qquad
 n_1=\mathbf E-\mathbf F,\qquad n_2=\mathbf E-\mathbf G.
\]
Their differentials are \((12,-4,-2)X\), so the highest cycles have basis
\(u=c+6n_2,\ v=n_1-2n_2\).
For \(w_b=f^{(b)}(x^2X)\), \(v_b=f^{(b)}(x^3X)\), take
\begin{align*}
 Z_2&=12w_3\otimes r_4-3w_2\otimes fr_4
                 +2w_1\otimes f^{(2)}r_4-3w_0\otimes f^{(3)}r_4,\\
 Z_3&=3v_3\otimes\rho_3-2v_2\otimes f\rho_3
                 +3v_1\otimes f^{(2)}\rho_3-12v_0\otimes f^{(3)}\rho_3 .
\end{align*}
After bracket-section elimination, the \(C_0,b_+\) head operators for
\(r_4\) are \(-C_1-3b_+\), \(3C_0\), and for \(\rho_3\) they are
\[
 -Q+b_-b_+/2,\qquad b_-C_0+C_1b_+/2-2\beta_+.
\]
The recursion \eqref{eq:boundary-divided-recursion} gives the following
two readings on each summand:
\[
\begin{array}{c|rr|rr}
 (a,b)&\mathbf A\text{ in }w_b r_a&\mathbf E\text{ in }w_b r_a&
       \mathbf A\text{ in }v_b\rho_a&\mathbf E\text{ in }v_b\rho_a\\\hline
 (0,3)&0&18&-6&96\\
 (1,2)&18&0&6&-156\\
 (2,1)&-18&-54&0&24\\
 (3,0)&0&36&0&36
\end{array}
\]
In a highest cycle \(au+bv\), these readings are \(a,6a-b\).
Multiplying the table by the coefficients in \(Z_2,Z_3\) gives
\[
 \partial_2Z_2=-90u-540v,\qquad
 \overline{\partial_2Z_3}=-30u-420v .
\]
Modulo seven these are \(u-v,5u\), removing both cycles.

For characteristic five, let \(\Omega=H^2+2H+4FE\) and
\(P=\Omega^2/9=I+N-N^2\) on weight one, where \(N=FE\).
The raw boundaries from \(y^2Y\otimes r_4\) and \(v_3\otimes\rho_3\) have
coordinates
\[
 s_2=(0,6,12,0,18,36,0),\qquad
 s_3=(-6,-6,-6,0,96,102,66).
\]
The raising and lowering formulas give
\[
\begin{array}{c|l}
 \xi&N\xi\\\hline
 \mathbf A&2\mathbf A+3\mathbf B+3\mathbf C\\
 \mathbf B&\mathbf A+4\mathbf B+\mathbf C+2\mathbf D\\
 \mathbf C&\mathbf A+\mathbf B+4\mathbf C+2\mathbf D\\
 \mathbf D&2\mathbf B+2\mathbf C+4\mathbf D\\
 \mathbf E,\mathbf F,\mathbf G&\mathbf E+\mathbf F+\mathbf G
\end{array}
\]
and hence \(Ps_2=2u+2v,\ Ps_3=4u+v\).
Their determinant is \(4\) modulo five.
The Casimir-three block of the target has three \(L(1)\) factors and
weights \(\pm1\). The relation \([E,F]=H\) makes it semisimple;
its cycle kernel is the two copies just removed.

At \((h,d,\zeta^\dagger)=(2,1,2)\) the invariant chains are
\begin{align*}
 \alpha_2={}&C_0\otimes y^2Y-C_1\otimes(y^2X+2xyY-tY)\\
            &+C_2\otimes(2xyX-tX+x^2Y)-C_3\otimes x^2X,\\
 \beta_2={}&b_+\otimes(y^2X-xyY)+b_-\otimes(x^2Y-xyX+tX),\\
 \gamma_2={}&t(b_+\otimes Y-b_-\otimes X).
\end{align*}
Their differential is \((-12,-5,-2)\) to \(\langle xY-yX\rangle\).
Thus the scalar cycle space is two-dimensional.
For \(w_j=f^{(j)}(x^3X)\), \(v_j=f^{(j)}(tx^2X)\), the invariant sources
\[
 \sum_{i=0}^4(12,-3,2,-3,12)_i\,w_{4-i}\otimes f^{(i)}r_4,\qquad
 \sum_{i=0}^3(3,-1,1,-3)_i\,v_{3-i}\otimes f^{(i)}\rho_3
\]
give the boundaries
\[
 -225\alpha_2+540\beta_2,\qquad
 15\alpha_2-60\beta_2+60\gamma_2 .
\]
Their summand readings in the ordered coordinates
\((C_0\otimes y^2Y,b_+\otimes y^2X,b_+\otimes tY)\) are
\[
\begin{array}{c|rrr|rrr}
 i&\multicolumn{3}{c|}{w_{4-i}\otimes f^{(i)}r_4}&
     \multicolumn{3}{c}{v_{3-i}\otimes f^{(i)}\rho_3}\\\hline
 0&-3&54&-36&3&-14&8\\
 1&39&-54&0&-6&52&-46\\
 2&-36&-54&108&0&-20&26\\
 3&0&54&-72&0&-18&12\\
 4&0&0&0&&&
\end{array}
\]
Modulo seven the first two coordinates of the two resulting boundaries
have determinant \(3\).
For characteristic five apply \(P=I+2\Omega-\Omega^3\) to the boundaries
\(s_2=\partial_2(y^3Y\otimes r_4)\),
\(s_3=\overline{\partial_2(y^4X\otimes\rho_3)}\).
Direct application of the displayed \(e,f\) formulas gives
\[
\begin{array}{c|rrr|rrr}
 &\multicolumn{3}{c|}{s_2}&\multicolumn{3}{c}{s_3}\\\hline
 I&2&4&4&4&2&4\\
 \Omega&3&3&4&3&3&0\\
 \Omega^2&1&1&0&2&4&0\\
 \Omega^3&3&1&2&4&2&0
\end{array}
\]
so \(Ps_2=4\beta_2,\ Ps_3=\alpha_2+\beta_2+4\gamma_2\).
These independent invariant boundaries fill the scalar cycle kernel.

At \((h,d,\zeta^\dagger)=(3,0,3)\), put
\begin{align*}
 \alpha_3={}&C_0\otimes y^3-3C_1\otimes(xy^2-yt)
          +3C_2\otimes(x^2y-xt)-C_3\otimes x^3,\\
 \beta_3={}&t(b_+\otimes y-b_-\otimes x).
\end{align*}
Their differential is \((-12,4)t\), leaving \(\alpha_3+3\beta_3\).
The same coefficient vector \((12,-3,2,-3,12)\), now with
\(w_j=f^{(j)}x^4\), has summand readings
\[
\begin{array}{c|rr}
 i&C_0\otimes y^3&b_+\otimes yt\\\hline
 0&-12&-144\\
 1&48&108\\
 2&-36&216\\
 3&0&-180\\
 4&0&0
\end{array}
\]
and gives \(-360(\alpha_3+3\beta_3)\), a unit multiple modulo seven.
For characteristic five take \(s=\partial_2(y^4\otimes r_4)\).
Its two readings under \(I,\Omega,\Omega^2,\Omega^3\) are
\((3,1),(0,3),(4,3),(4,0)\).
Thus \((I+2\Omega-\Omega^3)s=4\alpha_3+2\beta_3
=4(\alpha_3+3\beta_3)\).
All polynomial operators used here preserve the image of the boundary.

The \((h,d)=(1,2)\) entry follows from the quadratic cocycle equations
in Proposition~\ref{prop:h1-classification} below.
\end{proof}

\section{Tensor cocycles in the first two rows}\label{sec:first-two}

Put \(A=\kk[x,y,z]/(x^p,y^p,z^p)\),
\(\theta=dz-y\,dx+x\,dy\), and \(H^*=\Omega_A^1/A\theta\).
Write \(X,Y\) for its frame \(dx,dy\), and put
\[
 T(d,\lambda)=\operatorname{Sym}^dH^*\otimes\theta^k,\qquad
 k=(\lambda-d)/2,\qquad \lambda=-\zeta^\dagger .
\]
Any integer representative of \(k\in\mathbb F_p\) may be used.\footnote{Changing the representative by \(p\) does not change the infinitesimal action.} The PBW distribution pairing, as in
Section~\ref{sec:terminal-proof}, identifies \(T(d,\lambda)\) with
\(M(d,\zeta^\dagger)^*\): the negative fields act freely at the
origin, and the fibre is \(L(d)\) with central character \(\lambda\).
Let
\[
 E=X\partial_Y,\qquad F=Y\partial_X,\qquad
 H=X\partial_X-Y\partial_Y,\qquad \rho(I)=\lambda .
\]
The action is \(L_f=D_f+\rho(a_f)\), where the fibre matrices are
\begin{align}
 a_{C_0}&=6xE,& a_{C_1}&=-2xH+2yE,\notag\\
 a_{C_2}&=-2yH-2xF,&a_{C_3}&=-6yF,\notag\\
 a_{b_+}&=x(I+H)+2yE,&a_{b_-}&=2xF+y(I-H),\label{eq:tensor-action}\\
 a_Q&=2zI+2xyH-2x^2F+2y^2E.\notag
\end{align}
For example these follow by taking the horizontal derivatives of
\(D_fx=-f_y+xf_z\) and \(D_fy=f_x+yf_z\).
More explicitly, set
\(\nabla_x=\partial_x+y\partial_z\) and
\(\nabla_y=\partial_y-x\partial_z\).
Then \(L_f(dx)=d(D_fx)\) modulo \(\theta\), and similarly for \(dy\);
the action on \(\theta^k\) contributes \(2kf_z\).
These formulas also fix the signs of all cocycles below.
Internal degree \(r=h-1\) means that the value at an input of
contact degree \(j\) has coefficient degree \(r+j\).

\begin{lemma}\label{lem:degree-one-generation}
For every $p\ge5$, the positive contact algebra $\mathfrak n$ is generated by
\[
 \Kalg_1=\langle C_0,C_1,C_2,C_3,b_+,b_-\rangle.
\]
\end{lemma}
\begin{proof}
Its generated algebra contains all of
\(\Kalg_2=R_4\oplus zR_2\oplus\kk Q\), since
\[
 {[b_+,b_-]}=Q,\quad [C_0,b_+]=-x^4,\quad
 ([C_0,b_-]-[C_1,b_+])/4=zx^2 .
\]
Horizontal lowering fills \(R_4,zR_2\), using factorials at most \(4!\).
Induct on contact degree, and within it on the exponent of \(z\).
For a horizontal monomial of degree at least five, use
\[
 {[x^4,x^{a-3}y^{b+1}]}=4(b+1)x^ay^b
\]
if \(a\ge3,b\le p-2\), or its \(y^4\) counterpart
with coefficient \(-4(a+1)\).
The sole remaining corner is obtained from
\[
 {[Q,x^{p-2}y^{p-2}]}=-12zx^{p-2}y^{p-2},\qquad
 {[x^3,zx^{p-4}y^{p-1}]}
       =-3zx^{p-2}y^{p-2}-x^{p-1}y^{p-1}.
\]
For interior horizontal degree \(0<n<2p-2\), the identity
\[
 {[z\varphi,z^{k-1}g]}
   =z^k\{\varphi,g\}_{\rm sp}+(n-2)z^{k-1}\varphi g
       \quad(\varphi\in R_2,\ g\in R_n)
\]
reduces to smaller \(z\)-exponent.
Here \(\{R_2,R_n\}_{\rm sp}=R_n\): the coefficient of the
\(xy\)-action is \(b-a\), and when \(a=b\) one can instead use
\(\{x^2,x^{a-1}y^{a+1}\}_{\rm sp}=2(a+1)x^ay^a\).
All input contact degrees in this step are smaller.

For horizontal degree zero,
\([Q,z^{k-1}]=2(k-3)z^k\), except at \(k=3\);
that case follows from \([zy,z^2x]=-z^3+z^2xy\).
For horizontal degree \(2p-2\), the corresponding \(Q\)-coefficient
is \(2(k-5)\).
For \(p\ge7\), the exception \(k=5\) is settled by
\[
 {[x^3,z^6x^{p-4}y^{p-1}]}
   =-3z^6x^{p-2}y^{p-2}-6z^5x^{p-1}y^{p-1}.
\]
The first term comes from a \(Q\)-bracket with
coefficient \(-2\).
Every coefficient divided by in these steps is a unit.
\end{proof}

\begin{proposition}\label{prop:h0-classification}
In the row \(h=0\), the only nonzero groups are
\[
 (d,\lambda)=(0,0):L(1)\oplus L(3),\qquad
 (d,\lambda)=(1,1),(3,1):L(0).
\]
\end{proposition}
\begin{proof}
A cocycle is determined by the constant fibre values
\(a_i=\chi(C_i)\), \(b_\pm=\chi(b_\pm)\).
There are no inner cocycles of coefficient degree \(-1\).
The \(y\)-coefficient of the relation \(r_4\) gives
\(-8Ea_0=0\).
Transporting this equation by the horizontal group says that
the cubic-value map takes
\((sx+ty)^3\) into the line \(\kk(sX+tY)^d\).

We record the elementary pure-power calculation.
For a linear map from \(\operatorname{Sym}^mV\) to
\(\operatorname{Sym}^dV\) with this property, let \(A(s,t)\)
be its \(X^d\)-coefficient on \((sx+ty)^m\).
Its \(Y^d\)-coefficient on \(s\ne0\) is \(A(s,t)t^d/s^d\).
Polynomiality forces \(s^d\mid A(s,t)\).
Thus the map is determined by a polynomial of degree \(m-d\);
it is zero for \(d>m\), and the space of maps is \(L(m-d)\)
when \(d\le m<p\).
Consequently the cubic data lie in \(L(3-d)\), \(d\le3\).
If these data are zero, the same relation and its group translates
give the pure-power condition with \(m=1\) on the natural values.
Its kernel is zero for \(d\ge2\), one-dimensional for \(d=1\),
and two-dimensional for \(d=0\).

The scalar quadratic relation gives
\[
 Fa_0-Ha_1-Ea_2=0.
\]
On the possible \(d=1,2\) simple cubic components take
\(a_0=X^d\), \(a_1=(d/3)X^{d-1}Y\),
\(a_2=d(d-1)X^{d-2}Y^2/6\).
The left side is
\(2d(3-d)X^{d-1}Y/3\ne0\), excluding both.
For \(d=3\), the remaining cubic line is a scalar group module,
and its natural values are zero, since
\(\operatorname{Hom}_{\mathrm{SL}_2}(L(1),L(3))=0\).
Normalize \(a_i=X^{3-i}Y^i\). The \(x\)-coefficient of \(r_4\)
is \(-6(1+k)\), so \(k=-1\) and \(\lambda=1\).

For \(d=1\), the only possible nonzero data are
\(a_i=0,b_+=X,b_-=Y\). Put \(\Delta=xY-yX\).
The cocycle identity at \(Q=[b_+,b_-]\) gives
\(\chi(Q)=2(k-1)\Delta\).
Using
\[
 L_{b_+}\Delta=\lambda x\Delta-zX,\qquad
 L_QX=2\lambda zX-2x\Delta
\]
in \([b_+,Q]=0\) gives
\[
 \bigl(2(k-1)\lambda+2\bigr)x\Delta-6kzX=0 .
\]
Thus \(k=0,\lambda=1\).
For \(d=0,\lambda=2k\ne0\), the \(x\)-coefficient
\(-6ka_0\) of \(r_4\) first kills the cubic values.
Now \(\chi(Q)=\lambda(xb_--yb_+)\), and \([b_+,Q]=0\) gives
\[
 \lambda(\lambda+1)(x^2b_--xyb_+)-3\lambda z b_+=0 .
\]
Together with the corresponding negative equation this kills
both natural values.

We next construct the restricted classes attaining these bounds.
At \(d=\lambda=0\), every linear map \(\Kalg_1\to\kk\),
extended by zero on higher degrees, is a restricted cocycle:
\(\mathfrak n\) kills constants, and its \(p\)-powers have
degree at least \(p\). This gives \(\Kalg_1^*=L(1)\oplus L(3)\).

For \(d=1,\lambda=1\), the cocycle \(c_H\) of
Lemma~\ref{lem:contact-multiplier} has generator values
\(a_i=0,b_+=X,b_-=Y\). It gives the required nonzero invariant class.

For \(d=3,\lambda=1\), let
\(H=\ker\theta\), with frame \(\nabla_x,\nabla_y\),
and \(L=\operatorname{Der}(A)/H\), with frame
\(\ell=\partial_z\bmod H\).
The flat partial connection
\(\nabla^0_{\nabla_i}\nabla_j=0\) has Lie derivative
\(B_f=D_f\cdot\nabla^0\).
The connection identities make \(B_f\) bilinear over \(A\);
as the coboundary of a vector in a restricted Hom module it is
a restricted cocycle.
Lower its output by the Levi bracket
\(\omega(\nabla_x,\nabla_y)=-2\ell\), and symmetrize all three
indices. This is an equivariant tensor map, giving the restricted
cocycle
\[
 \chi(f)_{ijk}
       =2\operatorname{Sym}_{ijk}(\nabla_i\nabla_j\nabla_k f)\ell
       \quad\hbox{in }\operatorname{Sym}^3H^*\otimes L .
\]
The fibre character is \(3-2=1\).
Its values on \(C_i\) are \(12(dx)^{3-i}(dy)^i\ell\),
and it has internal degree \(-1\), so it is noninner.
Its pairing with the invariant cycle
\[
 \sum_{i=0}^3(-1)^i\binom3i C_i\otimes X^iY^{3-i}
\]
is \(48\): the fibre pairing of the corresponding terms is
\((-1)^i/\binom3i\). This is a unit for both primes.
The dimensions \(6,1,1\) attain the preceding upper bounds.
\end{proof}

\begin{proposition}\label{prop:h1-classification}
In the row \(h=1\), the only nonzero groups are
\[
 (d,\lambda)=(0,0):L(0),\qquad
 (d,\lambda)=(p-2,0):V^{(1)} .
\]
\end{proposition}
\begin{proof}
The generator values have coefficient degree one.
Write \(\chi(C_0)=xA+yA'\).
The \(z\)-coefficient of \(r_4\) is \(-3A'\), so \(A'=0\).
Transport by the horizontal group shows that
\(\chi((sx+ty)^3)\) is divisible by \(sx+ty\).
There are consequently unique fibre vectors \(A,B,C\) such that
\[
 \chi(C_0)=xA,\quad \chi(C_1)=(yA+2xB)/3,\quad
 \chi(C_2)=(xC+2yB)/3,\quad \chi(C_3)=yC .
\]
Define \(\eta(E)=A/6,\eta(H)=-B/3,\eta(F)=-C/6\);
then \(\chi(C_i)=\eta(a_{C_i})\).
For
\(\Theta(s,t)=\eta([s,t])-s\eta(t)+t\eta(s)\),
the scalar quadratic relation has value
\[
 12x^2\Theta(H,F)+12y^2\Theta(H,E)+24xy\Theta(E,F).
\]
It follows that \(\eta\) is an \(\mathfrak{sl}_2\)-cocycle.
Subtracting a constant inner tensor cocycle subtracts exactly
an inner \(\mathfrak{sl}_2\)-cocycle.

Here is its weight calculation.
Set \(\eta(E)=a,\eta(F)=b,\eta(H)=c\).
The three equations are
\[
 2a=Ha-Ec,\qquad -2b=Hb-Fc,\qquad c=Eb-Fa.
\]
In an integer weight \(\mu\) not divisible by \(p\), the cocycle
is inner, from \(c/\mu\).
The remaining weights are \(0,\pm p\), since
\(|\mu|\le d+2\le p+1\).
For even \(d\ge2\), weight zero has one-dimensional value spaces
at weights \(2,-2,0\); \(Ec=0\) forces \(c=0\), and
\(Eb=Fa\) leaves precisely the inner line.
For \(d=0\) perfectness gives zero.
For odd \(d\), only \(d=p-2,\mu=\pm p\) survive.
The positive representative is
\(\eta(F)=X^{p-2},\eta(E)=\eta(H)=0\);
the negative one is its group conjugate. There are no inner
vectors of these weights. Thus
\[
 H^1(\mathfrak{sl}_2,L(d))=
 \begin{cases}V^{(1)},&d=p-2,\\0,&d\ne p-2.\end{cases}
\]

We compute the kernel of the cubic-value map.
When \(\chi(C_i)=0\), write
\(\chi(b_+)=xP+yQ_0,\ \chi(b_-)=xR+yS\).
The \(r_4\) orbit gives
\begin{align*}
 Q_0+2EP&=0,&EQ_0&=0,&FR&=0,&R+2FS&=0,\\
 S-P+2ER-2HP&=0,&
 (1+2H)R+2FP&=0,\\
 2ES+(1-2H)Q_0&=0,&
 S-P-2HS-2FQ_0&=0 .
\end{align*}
These are the coefficients of \(s^3,s^2t,st^2,t^3\)
in the \(x,y\)-components of
\[
 (sx+ty)(s^2Q_0+st(S-P)-t^2R)
 +2(s^2E-stH-t^2F)
       \bigl(x(sP+tR)+y(sQ_0+tS)\bigr)=0 .
\]
If \(d\ge3\), the first two equations place \(P\) in
\(\langle X^d,X^{d-1}Y\rangle\), and the analogous equations place
\(S\) in \(\langle Y^d,XY^{d-1}\rangle\).
The supports in \(2ES+(1-2H)Q_0=0\) are disjoint:
they have \(Y\)-degrees \(d-2,d-1\) and zero.
Hence \(ES=0\), which forces \(S=R=0\);
then \(FP=0\) forces \(P=Q_0=0\).
For \(d=2\), write
\(P=p_0X^2+p_1XY,\ S=s_0Y^2+s_1XY\).
The equations give \(s_0=p_0=0,\ s_1=-3p_1,\ p_1=-3s_1\).
Since \(8\ne0\), all values vanish.

For \(d=1\) the kernel consists of
\[
 P=p_0X+p_1Y,\quad Q_0=-2p_1X,\quad
 R=2p_0Y,\quad S=-p_0X-p_1Y .
\]
It is a simple group module \(L(1)\), so a nonzero cocycle
kernel would include \(p_0=1,p_1=0\).
For this representative, putting \(\Delta=xY-yX\) gives
\(\chi(Q)=4kx\Delta\). The \(z\)-coefficient of \([b_+,Q]=0\)
is \(-4\lambda\), forcing \(\lambda=0,k=-1/2\).
The relation \(\rho_3\) then has value
\[
 L_{C_0}\chi(Q)+2L_{b_+}\chi(\beta_+)
                  -2L_{\beta_+}\chi(b_+)
      =(-6+18-6)x^3X=6x^3X\ne0,
\]
where \(\chi(\beta_+)=3x^2X\).
Thus no such cocycle exists.

For \(d=p-2\), the cubic-value map is now injective on cohomology,
with target \(V^{(1)}\).
A nonzero group submodule must be all of that simple module.
Normalize its weight-\(p\) representative as
\[
 \chi(C_0)=\chi(C_1)=0,\quad
 \chi(C_2)=-2xu_0,\quad\chi(C_3)=-6yu_0,\quad
 \chi(b_+)=0,\quad\chi(b_-)=a xu_0 .
\]
There is no constant inner vector of weight \(p\).
The \(xyu_1,y^2u_0\) coefficients of the lowest \(r_4\) relation
are \(18d(2-a)\) and \(9(4k-2-a)\).
They force \(a=2,k=1,\lambda=0\).

These two classes are restricted by the following construction.
In the divided-power module \(\Gamma^pH^*\) put
\(v_i=(dx)^{[p-i]}(dy)^{[i]}\).
Then
\[
 Ev_i=(p-i+1)v_{i-1},\qquad Fv_i=(i+1)v_{i+1}.
\]
The span \(\mathcal R\) of \(v_1,\ldots,v_{p-1}\) has basis
\(u_j=(-1)^jv_{j+1}\), with
\(Eu_j=ju_{j-1}\), \(Fu_j=(p-2-j)u_{j+1}\).
It is \(T(p-2,0)\), and there is a restricted module sequence
\[
 0\longrightarrow\mathcal R\longrightarrow\Gamma^pH^*
       \longrightarrow A\otimes V^{(1)}\longrightarrow0 .
\]
The constant quotient vectors have connecting cocycles
\[
 \chi_+(f)=-\nabla_y^2f\,u_0,\qquad
 \chi_-(f)=-\nabla_x^2f\,u_{p-2}.
\]
Restrictedness follows because divided powers are the symmetric-group
invariants in the tensor power of a restricted module.
If a nonzero constant quotient vector had an invariant lift,
its constant lift would be killed by both \(E,F\), since
\(C_0,C_3\) act by \(6xE,-6yF\).
But
\(\ker E=\langle v_0,v_1\rangle\) and
\(\ker F=\langle v_{p-1},v_p\rangle\) have zero intersection.
The connecting map is therefore injective on \(V^{(1)}\),
attaining the upper bound.

Finally, for \(d=0\) the kernel equations give
\(Q_0=R=0,P=S=D\).
For \(\lambda\ne0\) this is the inner cocycle of \(D/\lambda\).
For \(\lambda=0\), these are the generator values of \(Dc\)
from Lemma~\ref{lem:contact-multiplier}. The class of \(c\) attains the one-dimensional bound.
All other labels have been eliminated by ordinary relations, so the
same upper bounds hold for restricted cohomology.
\end{proof}

\section{Highest representatives and endpoint coordinates}\label{sec:highest-data}

We give the representative conventions used for the remaining rows.
For highest-vector arguments we may extend scalars to an algebraic
closure: the finite cochain complexes commute with field extension,
and vanishing descends by faithful scalar extension. The displayed
nonzero classes and their module maps are defined over the prime field.
All weights in this section are integer weights of the horizontal
torus. Set \(P_0=p-1\).

\begin{lemma}\label{lem:highest-lift}
Every highest cohomology class has a representative fixed by the
upper root group. The same conclusion holds for the quotient by
inner cochains of any group-stable space of generator data containing
the inner image.
\end{lemma}
\begin{proof}
For \(0\le s,d<p\), let \(E=x\partial_y+X\partial_Y\) on
\(L(s)\otimes L(d)\). Its kernel has basis
\[
 \eta_j=(xY-yX)^j x^{s-j}X^{d-j},
                  \qquad 0\le j\le\min(s,d).
\]
These are independent, nonzero, and killed by \(E\).
To see completeness, rescale the two monomial chains to identify
\(E\) with multiplication by \(u+v\) on
\(\kk[u,v]/(u^{s+1},v^{d+1})\).
Its cokernel has dimension \(\min(s,d)+1\), by setting \(v=-u\);
the kernel has the same dimension.
The weights \(s+d-2j\) of the kernel are nonnegative.
The invariant pairing on each simple factor is
\[
 B(x^{s-i}y^i,x^iy^{s-i})=(-1)^i/\binom si.
\]
It is nondegenerate, and \(B(Ea,b)=-B(a,Eb)\).
Hence the cokernel at weight \(\nu\) is dual to the kernel
at weight \(-\nu\). In particular \(E\) is onto every positive
weight space.

The coefficient space in degree \(r\) is
\[
 T_r(d,\lambda)=
 \bigoplus_{\substack{0\le c<p\\0\le r-2c\le2p-2}}
     z^cR_{r-2c}\otimes L(d).
\]
Here \(R_a\cong L(\min(a,2p-2-a))\), so the preceding
surjectivity holds on each summand.
Let \([\chi]\) be a highest class of weight \(w\ge0\).
If \(E\chi=\delta f\), choose \(g\) of weight \(w\)
with \(Eg=f\); then \(E(\chi-\delta g)=0\).
Divided root operators of order less than \(p\) are \(E^n/n!\)
and kill this representative.
Operators of order \(n\ge p\) take it into an inner cochain
of weight \(w+2n\ge2p\). The coefficient space has weights at
most \(p-1+d\le2p-2\), so these inner spaces are zero.
The representative is therefore group-highest.
The same argument applies to restricted cocycles.
\end{proof}

Here are explicit coordinates for these representatives.
Put \(m=\min(s,d)\), \(w_j=s+d-2j\), and
\(t_0=\max(0,s+d+2-p)\).
For \(j<t_0\), let
\[
 b_j=p-1+j-s,\qquad
 \gamma_j=y^sX^{d-b_j}Y^{b_j}.
\]
Expansion of the commuting summands of \(E\), using
\(\binom{p-1}{a}=(-1)^a\), gives
\begin{equation}\label{eq:long-endpoint-chain}
 E^{p-1}\gamma_j
       =(-1)^s s!b_j!\eta_j/j!\ne0,\qquad E^p\gamma_j=0.
\end{equation}
All factorial indices lie below \(p\).
For \(j\ge t_0\), one has \(0\le w_j<p\), and
\begin{equation}\label{eq:short-endpoint-chain}
 \pi_j=F^{w_j}\eta_j
      =w_j!(xY-yX)^j y^{s-j}Y^{d-j},\qquad
 E^{w_j}\pi_j=(w_j!)^2\eta_j .
\end{equation}
The latter follows inductively from
\(EF^a\eta_j=a(w_j-a+1)F^{a-1}\eta_j\).

These chains give a basis of the tensor product.
Their endpoints are independent; a kernel of the direct-sum map
would, by nilpotence of \(E\), contain a nonzero endpoint vector,
which is impossible.
Their total dimension, assuming \(s\le d\), is
\[
 pt_0+(s+1-t_0)(d+1-t_0)=(s+1)(d+1).
\]
This uses either \(t_0=0\) or \(t_0=s+d+2-p\).
Thus the chain lengths are \(p\) for \(j<t_0\) and \(w_j+1\)
otherwise.

For a function layer \(R_a\) above the horizontal truncation,
the concrete isomorphism needed in these formulas is
\begin{equation}\label{eq:truncated-horizontal-identification}
 x^{s-b}y^b\longmapsto
 \frac{\binom{p-1}{b}}{\binom{s}{b}}\,
           x^{p-1-b}y^{a-p+1+b},
 \qquad s=2p-2-a,\quad p\le a\le2p-2 .
\end{equation}
It intertwines \(E,F\) by adjacent binomial identities.
Below \(p\) use the identity map.
Apply this map to the function factor, leaving the fibre unchanged.

For a head module with basis \(H_0,\ldots,H_n\),
\(EH_i=iH_{i-1}\), take \(n=3\) for \(C_i\) and \(n=1\)
for \(b_+,b_-\). If a chain has bottom vector \(v_j\)
and length \(\ell_j\), then
\[
 v=E^{\ell_j-1-a}v_j,\qquad
 0\le a\le\min(n,\ell_j-1)
\]
runs through a basis of \(\ker E^{n+1}\).
The corresponding complete Lie-highest data are
\begin{equation}\label{eq:highest-endpoint-recovery}
 \chi_v(H_i)=\frac{i!}{n!}E^{n-i}v .
\end{equation}
The cocycle weight is \(w_j-2a+n\); retain its nonnegative values.

The group-highest condition has three possible edge corrections.
On monomials the target divided operator is
\[
 E_{\rm tar}^{(r)}(x^ay^bX^AY^B)
   =\sum_{u+v=r}\binom bu\binom Bv
        x^{a+u}y^{b-u}X^{A+v}Y^{B-v},
\]
with truncated monomials set to zero.
The \(p\)-th root coefficient on a cochain is
\[
 (E^{(p)}\chi)(H_i)=
    \sum_{a=0}^i(-1)^a\binom ia
              E_{\rm tar}^{(p-a)}\chi(H_{i-a}).
\]
For \(n=1\) the maximum data weight is \(2p-1\), so this
coefficient is zero on every nonnegative highest weight.
For \(n=3\), the maximum is \(2p+1\); only weights zero and one
require attention.
They occur at weight zero for
\((s,d)=(p-1,p-2),(p-2,p-1)\), and at weight one for
\((s,d)=(p-1,p-1)\).

In these cases write \(\chi_a\) for the chain datum with
\(j=p-a\) in \eqref{eq:highest-endpoint-recovery}.
The target weight after applying \(E^{(p)}\) is maximal and
one-dimensional; call its vector \(\chi_{\rm top}\).
The binomial expansion and
\(\binom{p+k}{p}=1\) for \(0\le k<p\) give
\[
 E^{(p)}\chi_a=(p-1-a)!\chi_{\rm top}.
\]
Indeed the mixed compositions
\(E_{\rm tar}^{(p-a)}E_{\rm tar}^{(a)}
 =\binom pa E_{\rm tar}^{(p)}\)
vanish for \(0<a<p\), and the remaining monomial is raised
through all of its \(y,Y\) powers.
It follows that the weight-zero basis is
\(\chi_2+3\chi_3\); the weight-one basis is
\[
 \chi_1+2\chi_2,\qquad \chi_1-6\chi_3.
\]
The respective factorial sums vanish modulo \(p\), and their
coefficients \(2,-6\) are units.
These give the kernels of the single nonzero root coefficient.
Finally \(E\chi=0\) makes the root polynomial a polynomial in
the \(p\)-th power of its parameter.
Its degree is at most \(2p+1\); after the \(p\)-th coefficient
vanishes, the group law kills the \(2p\)-th coefficient too,
since the coefficient of \(t^pu^p\) is twice that coefficient.
This proves completeness of the group-highest coordinates.

\subsection{The six coefficient equations}

For later coefficient readings write
\(\mathscr D(f,g)=L_f\chi(g)-L_g\chi(f)\).
The five \(L(4)\) relations and the scalar relation are
\begin{align}
 \mathcal R_0&=\mathscr D(C_0,C_1)+3\mathscr D(C_0,b_+),\notag\\
 \mathcal R_1&=2\mathscr D(C_0,C_2)+9\mathscr D(C_1,b_+)
                                      +3\mathscr D(C_0,b_-),\notag\\
 \mathcal R_2&=3\mathscr D(C_1,C_2)+\mathscr D(C_0,C_3)
                 +9\mathscr D(C_2,b_+)+9\mathscr D(C_1,b_-),\notag\\
 \mathcal R_3&=2\mathscr D(C_1,C_3)+3\mathscr D(C_3,b_+)
                                      +9\mathscr D(C_2,b_-),\notag\\
 \mathcal R_4&=\mathscr D(C_2,C_3)+3\mathscr D(C_3,b_-),\label{eq:six-quadratic-tests}\\
 \mathcal S&=\mathscr D(C_0,C_3)-3\mathscr D(C_1,C_2).\notag
\end{align}
Each underlying bracket sum is zero by
\eqref{eq:contact-bracket}.
The coefficient operators in these tests, besides the fibre
operators in \eqref{eq:tensor-action}, are
\begin{align*}
 L_{C_0}&=3x^2\partial_y-x^3\partial_z+6xE,\\
 L_{C_1}&=-x^2\partial_x+2xy\partial_y-x^2y\partial_z-2xH+2yE,\\
 L_{C_2}&=-2xy\partial_x+y^2\partial_y-xy^2\partial_z-2yH-2xF,\\
 L_{C_3}&=-3y^2\partial_x-y^3\partial_z-6yF,\\
 L_{b_+}&=x^2\partial_x+(z+xy)\partial_y+xz\partial_z
                                        +\lambda x+xH+2yE,\\
 L_{b_-}&=(-z+xy)\partial_x+y^2\partial_y+yz\partial_z
                                        +\lambda y-yH+2xF.
\end{align*}

\subsection{The \texorpdfstring{$Q$}{Q}-normal form}\label{sec:Q-normal}

Write \(R_a=\kk[x,y]/(x^p,y^p)_a\), setting \(R_a=0\) outside
\(0\le a\le2p-2\).
On the fibre \(L(d)\) put
\[
 \Delta=xY-yX,\qquad
 D_0=x\partial_X+y\partial_Y,\qquad
 \mathcal A=\Delta D_0 .
\]
The \(Q\)-operator in \eqref{eq:tensor-action} is
\begin{equation}\label{eq:Q-full-factorisation}
 L_Q=2z(x\partial_x+y\partial_y+z\partial_z+\lambda)-2\mathcal A .
\end{equation}
Indeed \(2xyH-2x^2F+2y^2E=-2\Delta D_0\).
Since \(D_0\Delta=0\), one has
\(\mathcal A^a=\Delta^aD_0^a\), and \(\mathcal A^a=0\) for \(a>d\).

Let \(1\le h\le2p-3\), and put \(\ell=\lfloor(h+1)/2\rfloor\).
For \(f=\sum_{j=0}^{\ell-1}z^jf_j\) of coefficient degree \(h-1\),
the coefficients of \(L_Qf\) are
\begin{equation}\label{eq:Q-triangular-coefficients}
 (L_Qf)_0=-2\mathcal A f_0,\qquad
 (L_Qf)_a=2(\lambda+h-a)f_{a-1}-2\mathcal A f_a
       \quad(1\le a\le\ell),\qquad f_\ell=0 .
\end{equation}
There is at most one zero diagonal, since \(\ell<p\).
In every nonzero row, descending elimination uses
\[
 f_{a-1}=\frac{g_a+2\mathcal A f_a}{2(\lambda+h-a)} .
\]
The horizontal degrees agree in this formula, and all \(z\)-exponents
remain below \(p\).

Let \(e\in\{0,\ldots,p-1\}\) represent \(\lambda+h\).
If \(e\notin\{1,\ldots,\ell\}\), subtracting a unique inner cochain
makes \(\chi(Q)\) horizontal, and \(L_Q\) on zero-cochains is injective.
If \(e=c\in\{1,\ldots,\ell\}\), one can instead arrange
\[
 \chi(Q)=q_0+z^cq_c,\qquad
 q_0\in R_{h+1}\otimes L(d),\quad
 q_c\in R_{h+1-2c}\otimes L(d).
\]
All remaining inner changes preserving this shape are generated by
\begin{equation}\label{eq:Q-residual-gauge}
 \Phi_c(u)=\sum_{a=0}^{c-1}
       \frac{z^{c-1-a}\mathcal A^a u}{a!},\qquad
 L_Q\Phi_c(u)=-\frac{2\mathcal A^cu}{(c-1)!}.
\end{equation}
The identity follows by cancellation of adjacent terms:
the \(z\)-multiplication coefficient at the \(a\)-th term is \(2a\).
Conversely the descending recurrence leaves only \(f_{c-1}=u\)
free, and reconstructs exactly this sum.
It follows that
\[
 \ker L_Q\cong\ker\mathcal A^c,\qquad
 \coker L_Q\cong
 \frac{R_{h+1}\otimes L(d)}
      {\mathcal A^c(R_{h+1-2c}\otimes L(d))}
       \oplus z^c(R_{h+1-2c}\otimes L(d)).
\]
The kernel isomorphism is induced by \(\Phi_c\).
These operations are horizontal-group equivariant.

\section{The remaining internal degrees}\label{sec:positive-fibres}

We use one relation matrix for all remaining fibres. Its columns are
the highest endpoint data of Section~\ref{sec:highest-data}; its rows
are the coefficient equations in contact degrees two and three.
The same matrix also treats the scalar fibre. We give its entries and
the finite reduction table in Appendix~\ref{app:comparison}.

For six generator values \(v=(v_0,\ldots,v_5)\), let
\(\mathcal Rv\) denote these quadratic and cubic defects.
Let \(\mathcal V_{h,d,w}\) be the space of group-highest generator
data of integer weight \(w\), with values of coefficient degree \(h\).
Write \(\mathcal J_{h,d,w}\) for the image of all weight-\(w\)
zero-cochains of degree \(h-1\), evaluated at the six generators.

\begin{lemma}[The finite comparison]\label{lem:finite-comparison}
Let \(p=5,7\), \(2\le h\le2p-3\), and \(0\le d<p\).
Modulo \(\mathcal V_{h,d,w}\cap\mathcal J_{h,d,w}\), the kernel
of \(\mathcal R\) on \(\mathcal V_{h,d,w}\) is zero except in the
following cases, where it has dimension one:
\begin{equation}\label{eq:cubic-comparison-table}
\begin{array}{c|c|c|c|l}
 h&d&\lambda&w&\text{generator data}\\ \hline
 p+1&0&0&p&\chi_x(f)=-x^{p-1}f_y\\
 2p-5&1&3&1&\text{order-two pole}\\
 2p-4&0&4&1&\chi(C_i)=0,\quad(\chi(b_+),\chi(b_-))=(H_y,-H_x)\\
 2p-4&1&1&0&\text{gradient pole}\\
 2p-3&0\le d<p&\text{all}&d&\text{top pole, }v=X^d.
\end{array}
\end{equation}
Here \(H=x^{p-1}y^{p-2}\). The order-two and gradient data are
specified below, and the top-pole data are
\eqref{eq:top-pole-data}. At \(h=p+1,d=\lambda=0\), the full
quadratic and cubic quotient, without the highest condition, has
dimension two, with weights \(p,-p\).
\end{lemma}

\begin{proof}
Use the monomial coefficient rules and bracket sections in
Appendix~\ref{app:comparison}. Table~\eqref{eq:cubic-comparison-table}
is its kernel table: the rank difference is
\eqref{eq:comparison-rank}, with the nonzero columns represented as
displayed. All omitted weights have rank difference zero.
The full-weight computation at the two Cartier labels is recorded
in \eqref{eq:cartier-full-comparison}.
We illustrate the coefficient elimination immediately below.
\end{proof}

\subsection{A quadratic kernel removed in degree three}

Take \(h=2,d=1\). The nonnegative highest weights are
\(6,4,2,0\). The inner sources \(xX,\Delta=xY-yX\) have readings
\(\delta(xX)(C_1)=-3x^2X\) and
\(\delta\Delta(C_0)=3x^2X\). After removing them, the complete
highest lists are
\begin{align*}
 w=6:&\quad v_3=tx^2X;\\
 w=4:&\quad v_2=ax^2X,\quad
 v_3=bx^2Y+(3a-b)xyX+ezX,\quad v_5=fx^2X;\\
 w=2:&\quad v_2=bx\Delta+dzX,\quad
 v_3=3by\Delta+3dzY,\\
 &\quad v_4=ux^2X,\quad v_5=qx^2Y+(u-q)xyX+szX;\\
 w=0:&\quad v_i=0\ (0\le i\le3),\quad
 v_4=ux\Delta+rzX,\quad v_5=uy\Delta+rzY.
\end{align*}
Unlisted values are zero. These lists have respectively
\(1,4,5,2\) parameters; they follow by listing the endpoint
monomials and imposing \(E^4v_3=E^2v_5=0\).

At weight six the reading \([zxX]\mathcal R_4=6t\) kills \(t\).
At weight four the three readings are
\[
\begin{array}{c|l}
 {[x^3X]}\mathcal S&21a+3b-e\\
 {[zxX]}\mathcal R_3&9a+3b-(10+3\lambda)e\\
 {[zyX]}\mathcal R_4&9a-3b-(2+3\lambda)e.
\end{array}
\]
Modulo seven these imply \(b=e=a=0\), and
\([x^2yX]\mathcal R_3=2f\) kills \(f\).
Modulo five they give \(b=3e,a=2e,(2-3\lambda)e=0\).
Only \(\lambda=4\) needs another pair of readings:
\([x^2yX]\mathcal R_3=f+3e\) and
\([x^3X]\mathcal R_2=e+4f\). They give \(f=2e,4e=0\).

At weight two the scalar readings \(24d,-9b+3d\) first kill
\(d,b\). The remaining three are
\[
 -18s,\qquad -63u+27q-9s,\qquad -27u+9q-3s.
\]
Modulo seven they successively kill \(s,q,u\); modulo five,
after \(s=0\), the matrix on \(u,q\) has determinant two.

At weight zero the quadratic equation is \(9u-3r=0\), hence
\(r=3u\). Put \(A_-=y^4,\beta_-=zy^2\). The cubic relation
\[
 C_3\wedge\beta_-+\tfrac14 C_2\wedge A_-
\]
has zero bracket, including at \(p=5\).
Its forced values are
\[
 \chi(A_-)=-L_{C_3}v_5=0,\qquad
 \chi(\beta_-)=(L_{C_2}v_5-L_{C_3}v_4)/4=2rzyY.
\]
Indeed the three terms before imposing \(r=3u\) are
\[
 L_{C_3}v_5=(3u-r)y^3Y,\quad
 L_{C_2}v_5=(3u-r)xy^2Y+2rzyY,\quad
 L_{C_3}v_4=(3u-r)y^3X-6rzyY.
\]
The cubic defect is \(-2ry^4Y\), so \(r=u=0\).
Thus the quadratic kernel in weight zero does not survive degree
three. This example also shows why the ordinary inner image and the
coefficient returns must be included in the comparison matrix.

\subsection{The two natural pole families}

\begin{lemma}\label{lem:order-two-pole}
In \(T(1,3)\), the following nonzero data at
\(h=2p-5\) do not extend to a cocycle:
\[
 \chi(C_i)=0,\quad
 \chi(b_+)=-9x^{p-1}y^{p-4}Y-3x^{p-2}y^{p-3}X,\quad
 \chi(b_-)=3x^{p-2}y^{p-3}Y+3x^{p-3}y^{p-2}X .
\]
\end{lemma}
\begin{proof}
Truncate \(x,y\), and retain the Laurent variable \(z\).
Put
\[
 \Omega=x^{p-1}y^{p-1},\quad H=x^{p-1}y^{p-2},\quad
 F=-\tfrac32(H_xX+H_yY),\quad\Theta=\tfrac92\Omega X,\quad
 P=F/z+\Theta/z^2 .
\]
Explicitly
\(F=3x^{p-1}y^{p-3}Y+\tfrac32x^{p-2}y^{p-2}X\).
Its coefficient degree is \(2p-6\).
The horizontal cubic action on \(F\) is zero, since
\(\{R_3,H\}_{\rm sp}=0\), and horizontal Lie differentiation
commutes with the ordinary differential.
The remaining cubic terms have horizontal degree beyond \(2p-2\).
Thus \(L_{C_i}P=0\).
The pole terms of the natural and \(Q\) actions cancel:
\[
\begin{array}{c|cc}
 &F/z&\Theta/z^2\\ \hline
 b_+&\frac92x^{p-1}y^{p-2}X/z&
                   -\frac92x^{p-1}y^{p-2}X/z\\
 b_-&-\frac92x^{p-2}y^{p-1}X/z&
                    \frac92x^{p-2}y^{p-1}X/z\\
 Q&9\Omega X/z&-9\Omega X/z .
\end{array}
\]
The polynomial natural values are those in the statement,
and \(L_QP=-4F\).
Brackets give polynomial values on degree two.
For an input of contact degree \(j\ge3\), a pole \(z^{-a}\)
would have horizontal degree
\(j+2p-6+2a\ge2p-1\), so it vanishes.
Equation~\eqref{eq:principal-part-defect} gives the defect
\[
 2(p-3)L_{z^p}P=12Fz^{p-2}+24\Theta z^{p-3}.
\]
The \(x^{p-1}y^{p-3}z^{p-2}Y\)-coefficient is \(36\ne0\).
We therefore exclude these data by Lemma~\ref{lem:principal-part}.
\end{proof}

Put \(d_Hf=f_xX+f_yY\).

\begin{lemma}\label{lem:gradient-pole}
In \(T(1,1)\), at \(h=2p-4\), the nonzero data
\[
 \chi(C_i)=0,\qquad \chi(b_+)=\tfrac12F_y,\qquad
 \chi(b_-)=-\tfrac12F_x,\qquad
 F=d_H(x^{p-1}y^{p-1})
\]
do not extend to a cocycle.
\end{lemma}
\begin{proof}
In the Laurent module take \(P=F/(2z)\), where
\[
 F=-x^{p-2}y^{p-1}X-x^{p-1}y^{p-2}Y .
\]
The horizontal cubic action on \(F\) is the differential of
\(\{C_i,x^{p-1}y^{p-1}\}_{\rm sp}=0\).
The remaining cubic products have excessive horizontal degree.
Thus \(L_{C_i}P=0\).
The pole numerator for \(b_+\) is
\(\bigl((2p-4+\lambda)x+xH+2yE\bigr)F\).
Since \(xF=-\Omega X\) and \((xH+2yE)F=-3\Omega X\),
it equals \((1-\lambda)\Omega X\), which is zero at \(\lambda=1\).
The \(b_-\) pole cancels in the same way.
The polynomial values are \(F_y/2,-F_x/2\).
Furthermore \(D_0F=-2\Omega\), so \(\mathcal A F=\Delta D_0F=0\);
the \(Q\)-action gives \(L_QP=-3F\).
Degree-two brackets have polynomial values.
For higher positive inputs, a pole would have horizontal degree at
least \(3+(2p-5)+2=2p>2p-2\).
Equation~\eqref{eq:principal-part-defect} gives the defect
\[
 2(p-3)L_{z^p}P=6Fz^{p-2}.
\]
The \(x^{p-2}y^{p-1}z^{p-2}X\)-coefficient is the unit \(-6\).
We exclude the stated data by Lemma~\ref{lem:principal-part}.
\end{proof}

\begin{proposition}\label{prop:scalar-columns}
For \(p=5,7\) and \(0\le h\le2p-3\), the scalar-fibre column is
\[
 H^1_{\mathrm{res}}(\mathfrak n,T(0,\lambda))_{h-1}
   \cong
 \begin{cases}
   L(1)\oplus L(3),&h=0,\ \lambda=0,\\
   L(0),&h=1,\ \lambda=0,\\
   V^{(1)},&h=p+1,\ \lambda=0,\\
   0,&\text{otherwise}.
 \end{cases}
\]
For every positive fibre degree \(1\le d<p\), the groups in
\(2\le h\le2p-3\) vanish at every density.
\end{proposition}
\begin{proof}
The first two rows were proved in
Propositions~\ref{prop:h0-classification} and~\ref{prop:h1-classification}.
For a remaining nonzero cohomology group choose a group-highest
class and lift it by Lemma~\ref{lem:highest-lift}. Its generator
values satisfy all quadratic and cubic equations. Subtracting an
inner cocycle puts these values in the corresponding line
of \eqref{eq:cubic-comparison-table}.

The scalar line at \(h=2p-4,\lambda=4\) is
\eqref{eq:scalar-pole-data}; its defect at \((Q,z^{p-1})\) is
\(12Hz^{p-2}\ne0\). The last line is excluded by
\eqref{eq:top-pole-defect}. Lemmas~\ref{lem:order-two-pole}
and~\ref{lem:gradient-pole} exclude the other natural lines.
In each case generation in degree one forces the higher values,
so subtracting an inner cocycle cannot change the obstruction.
There is consequently no nonzero class at any of these labels.

At \(h=p+1,d=\lambda=0\), the full quotient in
\eqref{eq:cartier-full-comparison} gives the upper bound two.
The two restricted Cartier cocycles constructed in
Theorem~\ref{thm:qtop} attain it and form \(V^{(1)}\).
These cases exhaust the comparison table.
\end{proof}

\begin{proof}[Proof of Theorem~\ref{thm:main}]
Theorem~\ref{thm:pbw-duality} gives
\eqref{eq:main-duality}, with \(r=h-1\).
The tensor realization \eqref{eq:tensor-action} identifies its
coefficient calculation with the tensor cocycles used above.
Propositions~\ref{prop:h0-classification} and~\ref{prop:h1-classification}
give the first five lines of \eqref{eq:corrected-main-table}.
Their restricted cocycles attain the ordinary upper bounds.
Proposition~\ref{prop:scalar-columns} gives the remaining rows,
including the Cartier row of Theorem~\ref{thm:qtop}.
All nonzero horizontal modules are self-dual, so finite duality
gives exactly the modules in \eqref{eq:corrected-main-table}.
Finally \(q=2p\) gives \(h=2p-3\), which for \(p=5,7\)
is neither zero, one, nor \(p+1\). This proves its asserted vanishing.
\end{proof}

\appendix
\section{The coefficient matrix for the finite comparison}\label{app:comparison}

This appendix specifies the matrix and the reduction rule for
\eqref{eq:cubic-comparison-table}. All parameters are in the finite
ranges \(p=5,7\), \(2\le h\le2p-3\), \(0\le d<p\) and
\(\lambda\in\mathbb F_p\). Highest columns have \(0\le w\le2p+1\);
for the full Cartier comparison we take \(|w|\le2p+1\).
All weights are integers.

\subsection{Monomial entries}

Order the basis of \(T_r(d,\lambda)\) by \((c,a,b,i)\), where
\[
 m_{a,b,c,i}=x^ay^bz^cX^{d-i}Y^i,\qquad
 a+b+2c=r,\quad0\le a,b,c<p,\quad0\le i\le d.
\]
Its weight is \(a-b+d-2i\); write \(T_r[w]\) for the corresponding
weight subspace. Any monomial outside these bounds is zero.
The input order is
\(e_0=C_0,e_1=C_1,e_2=C_2,e_3=C_3,e_4=b_+,e_5=b_-\).
For arbitrary monomial inputs the bracket is
\begin{align}
 [x^ay^bz^c,x^Ay^Bz^C]
 ={}&(aB-bA)x^{a+A-1}y^{b+B-1}z^{c+C}\notag\\
 &+\bigl((2-a-b)C+c(A+B-2)\bigr)
                  x^{a+A}y^{b+B}z^{c+C-1}.
 \label{eq:comparison-monomial-bracket}
\end{align}
For all coefficient actions, put
\[
 P_x=-f_y+xf_z,\quad P_y=f_x+yf_z,\quad
 P_z=2f-xf_x-yf_y,
 \qquad \nabla_x=\partial_x+y\partial_z,\quad
 \nabla_y=\partial_y-x\partial_z.
\]
The required entry formula is
\begin{align}
 L_f={}&P_x\partial_x+P_y\partial_y+P_z\partial_z
       +(\lambda-d)f_z\notag\\
 &+(X\nabla_xP_x+Y\nabla_yP_x)\partial_X
  +(X\nabla_xP_y+Y\nabla_yP_y)\partial_Y.
 \label{eq:comparison-action}
\end{align}
Thus each entry is obtained by differentiation and multiplication
of the displayed monomials, followed by truncation. For example
\(L_{C_0}m\) has three terms, with coefficients
\(3b,-c,6i\), and exponent changes
\((2,-1,0,0),(3,0,-1,0),(1,0,0,-1)\).
This agrees with the six operators in
Section~\ref{sec:highest-data}.

\subsection{Highest columns and inner columns}

At weight \(w\ge0\), begin with independent endpoints
\[
 P\in T_h[w-3],\qquad Q_-\in T_h[w-1],
 \qquad v_i=\frac{i!}{6}E^{3-i}P\ (0\le i\le3),
 \quad v_4=EQ_-,\quad v_5=Q_-.
\]
Here \(E=x\partial_y+X\partial_Y\) acts on coefficients.
Impose \(Ev_i=iv_{i-1}\) on each input string, with the initial
value zero, and impose its \(p\)-th root equation
\[
 \sum_{a=0}^i(-1)^a\binom ia E_{\rm tar}^{(p-a)}v_{i-a}=0.
\]
For the natural string replace \(v_i\) by \(v_{4+i}\),
\(0\le i\le1\). The monomial divided action is given in
Section~\ref{sec:highest-data}. These are all group-highest equations.
Let \(V\) be the matrix whose columns are a basis of their kernel,
expressed as six generator values. There are at most thirty-two
such columns in the stated ranges.

The inner matrix \(J\) has, for each monomial
\(m\in T_{h-1}[w]\), the column
\[
 (L_{e_0}m,\ldots,L_{e_5}m).
\]
In a common row basis of six coefficient spaces,
the dimension of their intersection with the highest space is
\begin{equation}\label{eq:comparison-inner-intersection}
 \dim(\im V\cap\im J)
     =\dim\mathcal V_{h,d,w}+\operatorname{rank}J
                                    -\operatorname{rank}[V\ J].
\end{equation}

\subsection{Quadratic and cubic rows}

Use the six equations \eqref{eq:six-quadratic-tests} for the
quadratic rows. The cubic rows can be specified without choosing
individual relations. Let \(\mathcal B_j\) be the lexicographically
ordered monomials \(x^ay^bz^c\), \(a+b+2c=j+2\), \(0\le a,b,c<p\).
The bracket matrices are
\[
 B_2=( [e_i,e_j])_{i<j},\qquad
 B_3=( [e_i,f])_{0\le i<6,\ f\in\mathcal B_2},
\]
with rows \(\mathcal B_2,\mathcal B_3\), respectively.
Order the columns lexicographically by their indicated pairs.
Their independent column lists, numbered from zero, are
\begin{equation}\label{eq:comparison-section-pivots}
\begin{array}{c|l}
 B_2&0,1,2,4,6,8,9,11,14\\
 B_3,\ p=5&0,1,2,3,4,10,11,19,37,39\\
 B_3,\ p=7&0,1,2,3,4,5,10,11,19,20,37,39.
\end{array}
\end{equation}
The corresponding square matrices \(A_j\) are invertible by
\eqref{eq:comparison-monomial-bracket}.
For a six-tuple \(v\) first form
\[
 d_{ij}=L_{e_i}v_j-L_{e_j}v_i.
\]
If \(I_2\) is the first list in
\eqref{eq:comparison-section-pivots}, define its forced
degree-two values by
\[
 t_f=\sum_{\alpha\in I_2}(A_2^{-1})_{\alpha,f}d_\alpha
                    \qquad(f\in\mathcal B_2).
\]
Next form
\[
 d_{i,f}=L_{e_i}t_f-L_fv_i,\qquad
 s_g=\sum_{\alpha\in I_3}(A_3^{-1})_{\alpha,g}d_\alpha
                    \qquad(g\in\mathcal B_3).
\]
There are fifty-four cubic equations
\begin{equation}\label{eq:comparison-cubic-rows}
 d_{i,f}-\sum_{g\in\mathcal B_3}(B_3)_{g,(i,f)}s_g=0.
\end{equation}
Several are identically zero. Let \(D\)
be the matrix of these equations and the six quadratic equations
on the highest columns \(V\).

An inner column satisfies all these equations by the module
identity. Since \(V\) has independent columns, the kernel dimension
after inner normalization is therefore
\begin{equation}\label{eq:comparison-rank}
 \nu_{p,h,d,\lambda,w}
 =\operatorname{rank}[V\ J]-\operatorname{rank}J
                                      -\operatorname{rank}D.
\end{equation}
This proves that the matrix measures precisely the quotient in
Lemma~\ref{lem:finite-comparison}.

\subsection{Reduction table}

The nonzero entries in these ranges are
\[
\begin{array}{c|c|c|c|c|c}
 p&h&d&\lambda&w&\nu\\ \hline
 5&5&1&3&1&1\\
 5&6&0&0&5&1\\
 5&6&0&4&1&1\\
 5&6&1&1&0&1\\
 5&7&0,1,2,3,4&0,1,2,3,4&d&1\\ \hline
 7&8&0&0&7&1\\
 7&9&1&3&1&1\\
 7&10&0&4&1&1\\
 7&10&1&1&0&1\\
 7&11&0,1,\ldots,6&0,1,\ldots,6&d&1.
\end{array}
\]
In the last line for each prime the \(d\) and \(\lambda\)
entries vary independently. All other entries are zero.
The remaining highest columns are the five representatives in
\eqref{eq:cubic-comparison-table}, with no linear combination
belonging to \(\im J\).

Here is a uniform rule for reading or checking any entry of the table.
Reduce coefficients modulo \(p\) after each operation. Solve the
highest equations by ascending endpoint columns, expressing pivot
coordinates in terms of the free coordinates. Form \(V,J,D\) as
above. In each matrix used in \eqref{eq:comparison-rank}, scan the
columns in order. For a new pivot choose the first unused row with
nonzero entry, divide that row by its pivot and subtract it from
the other rows. Count one for each pivot. Empty matrices have rank
zero. The elimination for \(h=2,d=1\) in the text
gives a complete example, including its cubic equation.

At the Cartier labels one may omit the highest equations and use
all six generator values of a fixed weight as independent columns.
If \(D_w^{\rm all}\) is the resulting relation matrix, the full
quotient has weight dimensions
\begin{equation}\label{eq:cartier-full-comparison}
 \dim\mathcal C_w
 =\dim\!\left(\bigoplus_{i=0}^5T_{p+1}[w+\operatorname{wt}e_i]\right)
       -\operatorname{rank}D_w^{\rm all}-\operatorname{rank}J_w
 =\begin{cases}1,&w=p,-p,\\0,&w\ne p,-p.\end{cases}
\end{equation}
For \(w=p,-p\), the column dimension is five. The inner matrix
\(J_w\) is zero, and \(D_w^{\rm all}\) has rank four.
For the other weights their two ranks sum to the column dimension.
These statements use the same coefficient and section rules.
The two columns are the contractions of
\(x^{p-1}dx,y^{p-1}dy\); their noninnerness was proved in
Section~\ref{sec:terminal-proof}.

\section*{Acknowledgements}

I thank J\"org Feldvoss for many helpful comments on this series of papers.

\end{document}